\documentclass[11pt]{amsart}
\usepackage{latexsym,graphicx}

\numberwithin{equation}{section}
\theoremstyle{plain}

\theoremstyle{remark}

\theoremstyle{definition}

\newcommand{\D}{{\mathcal D}}
\newcommand{\E}{\mathcal E}

\newcommand{\G}{{\mathcal G}}

\newcommand{\K}{{\mathcal K}}
\renewcommand{\L}{{\mathcal L}}
\newcommand{\M}{{\mathcal M}}
\newcommand{\N}{\mathbb N}

\newcommand{\V}{{\mathcal V}}

\newcommand{\dist}{\operatorname{dist}}

\newcommand{\fp}{\operatorname{FP}}

\newcommand{\Int}{\operatorname{Int}}

\renewcommand{\span}{\operatorname{span}}

\newcommand{\supp}{\operatorname{Supp}}

\def\half{{1 \over 2}}

\newcommand{\oa}{\overrightarrow}

\newcommand{\ol}{\overline}

\def\XXint#1#2#3{{\setbox0=\hbox{$#1{#2#3}{\int}$}
      \vcenter{\hbox{$#2#3$}}\kern-.5\wd0}}

\begin{document}

\def\cal{\mathcal}

\font\tpt=cmr10 at 12 pt
\font\fpt=cmr10 at 14 pt

\font \fr = eufm10




\overfullrule=0in

\def\boxit#1{\hbox{\vrule
 \vtop{%
  \vbox{\hrule\kern 2pt %
     \hbox{\kern 2pt #1\kern 2pt}}%
   \kern 2pt \hrule }%
  \vrule}}

  \def\harr#1#2{\ \smash{\mathop{\hbox to .3in{\rightarrowfill}}\limits^{\scriptstyle#1}_{\scriptstyle#2}}\ }

\def\AA{1}
\def\BB{2}
\def\CC{3}
\def\DD{4}
\def\EE{5}
\def\FF{6}
\def\GGG{7}
\def\HH{8}
\def\II{9}
\def\JJ{10}
\def\KK{11}
\def\LL{12}
\def\MM{13}

 \def\GG{{{\bf G} \!\!\!\! {\rm l}}\ }

\def\GL{{\rm GL}}

\def\bll{I \!\! L}

\def\IFF{\qquad\iff\qquad}
\def\bra#1#2{\langle #1, #2\rangle}
\def\bbf{{\bf F}}
\def\bbj{{\bf J}}
\def\Jtn{{\bbj}^2_n}  \def\JtN{{\bbj}^2_N}  \def\JoN{{\bbj}^1_N}
\def\jt{j^2}
\def\jtx{\jt_x}
\def\Jt{J^2}
\def\Jtx{\Jt_x}
\def\bpp{{\bf P}^+}
\def\bpt{{\wt{\bf P}}}
\def\fsh{$F$-subharmonic }
\def\mo{monotonicity }
\def\jet{(r,p,A)}
\def\ss{\subset}
\def\sse{\subseteq}
\def\half{\hbox{${1\over 2}$}}
\def\smfrac#1#2{\hbox{${#1\over #2}$}}
\def\oa#1{\overrightarrow #1}
\def\dim{{\rm dim}}
\def\dist{{\rm dist}}
\def\codim{{\rm codim}}
\def\deg{{\rm deg}}
\def\rank{{\rm rank}}
\def\log{{\rm log}}
\def\Hess{{\rm Hess}}
\def\Hessyp{{\rm Hess}_{\rm SYP}}
\def\trace{{\rm trace}}
\def\tr{{\rm tr}}
\def\max{{\rm max}}
\def\min{{\rm min}}
\def\span{{\rm span\,}}
\def\Hom{{\rm Hom\,}}
\def\det{{\rm det}}
\def\End{{\rm End}}
\def\Sym{{\rm Sym}^2}
\def\diag{{\rm diag}}
\def\pt{{\rm pt}}
\def\Spec{{\rm Spec}}
\def\pr{{\rm pr}}
\def\Id{{\rm Id}}
\def\Grass{{\rm Grass}}
\def\Herm#1{{\rm Herm}_{#1}(V)}
\def\arr{\longrightarrow}
\def\supp{{\rm supp}}
\def\Link{{\rm Link}}
\def\Wind{{\rm Wind}}
\def\Div{{\rm Div}}
\def\vol{{\rm vol}}
\def\foral{\qquad {\rm for\ all\ \ }}
\def\fpsh{{\cal PSH}(X,\f)}
\def\Core{{\rm Core}}
\def\dis{f_M}
\def\Re{{\rm Re}}
\def\rn{\bbr^n}
\def\pp{\cp^+}
\def\plp{\cp_+}
\def\Int{{\rm Int}}
\def\cix{C^{\infty}(X)}
\def\Gr#1{G(#1,\rn)}
\def\Symn{{\Sym(\rn)}}
\def\SymN{{\Sym(\bbr^N)}}
\def\Gpn{G(p,\rn)}
\def\fd{{\rm free-dim}}
\def\SA{{\rm SA}}
 \def\cd{{\cal C}}
 \def\cdt{{\widetilde \cd}}
 \def\cm{{\cal M}}
 \def\cmt{{\widetilde \cm}}

\def\Theorem#1{\medskip\noindent {\bf THEOREM \bf #1.}}
\def\Prop#1{\medskip\noindent {\bf Proposition #1.}}
\def\Cor#1{\medskip\noindent {\bf Corollary #1.}}
\def\Lemma#1{\medskip\noindent {\bf Lemma #1.}}
\def\Remark#1{\medskip\noindent {\bf Remark #1.}}
\def\Note#1{\medskip\noindent {\bf Note #1.}}
\def\Def#1{\medskip\noindent {\bf Definition #1.}}
\def\Claim#1{\medskip\noindent {\bf Claim #1.}}
\def\Conj#1{\medskip\noindent {\bf Conjecture \bf    #1.}}
\def\Ex#1{\medskip\noindent {\bf Example \bf    #1.}}
\def\Qu#1{\medskip\noindent {\bf Question \bf    #1.}}
\def\Exercise#1{\medskip\noindent {\bf Exercise \bf    #1.}}

\def\HoQu#1{ {\AAA T\BBB HE\ \AAA H\BBB ODGE\ \AAA Q\BBB UESTION \bf    #1.}}

\def\pf{\medskip\noindent {\bf Proof.}\ }
\def\qed{\hfill  $\vrule width5pt height5pt depth0pt$}
\def\equdef{\buildrel {\rm def} \over  =}
\def\qedqed{\hfill  $\vrule width5pt height5pt depth0pt$ $\vrule width5pt height5pt depth0pt$}
\def\mathqed{  \vrule width5pt height5pt depth0pt}

\def\V{W}

\def\df{d^{\phi}}
\def\hk{\_{\rm l}\,}
\def\n{\nabla}
\def\w{\wedge}

\def\cu{{\cal U}}   \def\cc{{\cal C}}   \def\cb{{\cal B}}  \def\cz{{\cal Z}}
\def\cv{{\cal V}}   \def\cp{{\cal P}}   \def\ca{{\cal A}}
\def\cw{{\cal W}}   \def\co{{\cal O}}
\def\ce{{\cal E}}   \def\ck{{\cal K}}
\def\ch{{\cal H}}   \def\cm{{\cal M}}
\def\cs{{\cal S}}   \def\cn{{\cal N}}
\def\cd{{\cal D}}
\def\cl{{\cal L}}
\def\cp{{\cal P}}
\def\cf{{\cal F}}
\def\ccr{{\cal  R}}

\def\gerG{{\fr{\hbox{g}}}}
\def\gerB{{\fr{\hbox{B}}}}
\def\gerR{{\fr{\hbox{R}}}}
\def\p#1{{\bf P}^{#1}}
\def\vf{\varphi}

\def\wt{\widetilde}
\def\wh{\widehat}

\def\and{\qquad {\rm and} \qquad}
\def\arr{\longrightarrow}
\def\ol{\overline}
\def\bbr{{\mathbb R}}\def\bbh{{\mathbb H}}\def\bbo{{\mathbb O}}
\def\bbc{{\mathbb C}}
\def\bbq{{\mathbb Q}}
\def\bbz{{\mathbb Z}}
\def\bbp{{\mathbb P}}
\def\bbd{{\mathbb D}}

\def\a{\alpha}
\def\b{\beta}
\def\d{\delta}
\def\e{\epsilon}
\def\f{\phi}
\def\g{\gamma}
\def\k{\kappa}
\def\l{\lambda}
\def\o{\omega}

\def\s{\sigma}
\def\x{\xi}
\def\z{\zeta}

\def\D{\Delta}
\def\L{\Lambda}
\def\G{\Gamma}
\def\O{\Omega}

\def\bd{\partial}
\def\bdf{\partial_{\f}}
\def\lag{Lagrangian}
\def\psh{plurisubharmonic }
\def\ph{pluriharmonic }
\def\pph{partially pluriharmonic }
\def\omp{$\omega$-plurisubharmonic \ }
\def\ffl{$\f$-flat}
\def\PH#1{\widehat {#1}}
\def\lloc{L^1_{\rm loc}}
\def\dbar{\ol{\partial}}
\def\lp{\Lambda_+(\f)}
\def\lpp{\Lambda^+(\f)}
\def\bo{\partial \Omega}
\def\Ob{\overline{\O}}
\def\fc{$\phi$-convex }
\def\PSH{{ \rm PSH}}
\def\SH{{\rm SH}}
\def\totr{ $\phi$-free }
\def\BM{\lambda}
\def\Der{D}
\def\CH{{\cal H}}
\def\RH{\overline{\ch}^\f }
\def\pconv{$p$-convex}
\def\MA{MA}
\def\lagpsh{Lagrangian plurisubharmonic}
\def\hermsk{{\rm Herm}_{\rm skew}}
\def\PSHl{\PSH_{\rm Lag}}
 \def\ppsh{$\pp$-plurisubharmonic}
\def\fp{$\pp$-plurisubharmonic }
\def\fh{$\pp$-pluriharmonic }
\def\Symn{\Sym(\rn)}
 \def\ci{C^{\infty}}
\def\USC{{\rm USC}}
\def\LSC{{\rm LSC}}
\def\fa{{\rm\ \  for\ all\ }}
\def\ppc{$\pp$-convex}
\def\cpt{\wt{\cp}}
\def\ft{\wt F}
\def\ob{\overline{\O}}
\def\Be{B_\e}
\def\K{{\rm K}}

\def\M{{\bf M}}
\def\N#1{C_{#1}}
\def\ds{Dirichlet set }
\def\dir{Dirichlet }
\def\Fa{{\oa F}}
\def\TR{{\cal T}}
 \def\ISO{{\rm ISO_p}}
 \def\Span{{\rm Span}}

\def\ALL{1}
\def\AV{2}
\def\BTA{3}
\def\BL{4}
\def\BRE{5}
\def\CNS{6}
\def\CP{7}
\def\CPW{8}
\def\CIL{9}
\def\CRA{10}
\def\DTT{11}
\def\DON{12}
\def\CG{13}
\def\DDD{14}
\def\DDR{15}
\def\GEO{16}
\def\HYP{17}
\def\BEL{18}
\def\SURVEY{19}
\def\AC{20}
\def\NOTES{21}
\def\AET{22}
\def\TANG{23}
\def\TANGG{24}
\def\LAG{25}
\def\SLE{26}
\def\JTY{27}
\def\KRY{28}
\def\PLI{29}
\def\RT{30}
\def\SLO{31}
\def\SPR{32}
\def\TRUU{33}
\def\TRU{34}
\def\TWC{35}
\def\TWCC{36}
\def\TWCCC{37}
\def\WAL{38}

\def\II{1}
\def\AA{2}
\def\AB{3}
\def\EE{4}
\def\BB{5}
\def\CC{6}
\def\CCC{7}
\def\DD{8}

\vskip .4in

\def\E{E}
\def\fpsi{{F_f(\psi)}}
\def\bL{{\bf \Lambda}}
\def\bdf{{\bf f}}
\def\UU{U}
\def\bbm{{\bf M}}

\font\headfont=cmr10 at 14 pt
\font\aufont=cmr10 at 11 pt

\title[THE INHOMOGENEOUS DIRICHLET PROBLEM]
{\headfont THE INHOMOGENEOUS DIRICHLET PROBLEM\\  
FOR GEOMETRIC SUBEQUATIONS ON \\  MANIFOLDS}

\title[THE INHOMOGENEOUS DIRICHLET PROBLEM]
{\headfont THE INHOMOGENEOUS DIRICHLET PROBLEM\\  
 for NATURAL  OPERATORS ON MANIFOLDS}

\date{\today}
\author{ F. Reese Harvey and H. Blaine Lawson, Jr.}
\thanks
{Partially supported by the NSF}

\maketitle

\centerline
{\sl   Dedicated with great admiration to Marcel Berger }

\begin{abstract}
We shall discuss the inhomogeneous Dirichlet problem for:
$$
f(x,u, Du, D^2u) \ =\ \psi(x)
$$
where $f$ is a ``natural'' differential operator, with a restricted domain $F$,
on a manifold $X$.  By ``natural''  we mean operators that arise intrinsically from
a given geometry on $X$.  An important point is that the equation need not be convex and can be highly degenerate.  Furthermore, the  inhomogeneous term can take values at the boundary of the restricted  domain $F$ of the operator $f$.

A simple example is the real Monge-Amp\`ere operator
$
\det\left( \Hess\, u  \right) \ =\ \psi(x)
$
on a riemannian manifold $X$, where $\Hess$ is the riemannian Hessian,  the  restricted domain is
 $F = \{\Hess \geq 0\}$, and $\psi$ is continuous with $\psi\geq0$.

A  main new tool is the idea of local jet-equivalence,
which gives rise to local weak comparison, and then to comparison under
a natural and necessary global assumption.

The main theorem applies to pairs $(F,f)$, which are locally jet-equivalent to a given constant coefficient pair $(\bbf, \bdf)$.
This covers a large family of geometric equations on manifolds: orthogonally invariant operators on a riemannian manifold,
G-invariant operators on manifolds with G-structure, operators on almost complex manifolds,
and operators, such as the Lagrangian Monge-Amp\`ere operator, on symplectic manifolds.
It also applies to all branches of these operators.
  Complete existence and uniqueness results are established
with existence requiring the same boundary assumptions as in the homogeneous case [\DDR].

We also have results where the inhomogeneous term $\psi$ is a delta function.
\end{abstract}


\vskip .1in

\vfill \eject


\
\vskip .5in

\centerline{\bf Table of Contents}  
\vskip.3in


 \hskip .5 in  \II.  Preliminary Discussion
 \medskip

 \hskip .5 in  \AA.  Statement of the Main Result
 \medskip

 \hskip .5 in  \AB.   Boundary Convexity
 \medskip

 \hskip .5 in  \EE.  Local Weak Comparison
 \medskip

 \hskip .5 in  \BB. Proof of the Main Theorem
 \medskip

 \hskip .5 in  \CC. Applications and Historical Remarks
 \medskip

 \hskip .5 in  \CCC. Fundamental Solutions
 \medskip

 \hskip .5 in  \DD. Two Important Questions

\vskip .3in

\ \qquad \qquad Appendix A.  Comparison for Constant Coefficient Operators

\ \qquad \qquad  Appendix B. Compatibility and Topological  Tameness
 are \\ . \qquad\qquad\qquad\qquad\qquad\qquad Necessary Conditions


\vfill\eject

\noindent {\headfont \II.    Preliminary Discussion.}     

The objective of this paper is to provide  a solution to the Dirichlet Problem 
with continuous boundary data and inhomogeneous term for a wide class of 
geometrically interesting equations on manifolds. 
The main points are that the equation need not be either  convex or invariant (as in [\CNS]), 
and it is allowed to be highly degenerate.  Complete existence and
uniqueness results are established.  
Comparison, and hence uniqueness, is proven under a mild strengthening
 of the standard weak ellipticity  assumption
on the operator $f$, which we call {\sl tameness}.
Existence requires the same boundary
assumption as in the homogeneous case  [\DDR].

The operators  considered here are those which are locally jet-equivalent
to a constant coefficient, or {\sl eucildean}, case.   The notion of jet-equivalence, introduced
in [\DDR], is very general.  It is not like transformations of coordinates; it almost never
takes the 2-jet of a function to a 2-jet of any function.  However, it is exactly
what is necessary for treating interesting geometric equations on manifolds.
The reader may want to look at the examples below, and in \S \CC.

The results here are a direct extension of the work in [\DDR]  in the following sense.
Here we assume that the differential operator and its domain $(F,f)$  are locally jet-equivalent
to a constant coefficient pair $(\bbf, \bdf)$.  Results in [\DDR] then establish the Dirichlet problem
for $f(J^2u)=0$ (where $J^2u$ is the 2-jet of $u$). In fact in [\DDR] there were no operators.
We simply replaced the pair $(F,f)$ with the subequation $F_f \equiv \{f\geq0\}$ and took a potential theory point-of-view.
Here we consider the problem
$$
f(J^2u) \ =\ \psi
$$
 where $\psi$ is an arbitrary continuous function  with values in the range of $f$ ($=\bdf(\bbf)$). 
  The principle work which reduces this to certain results in [\DDR]
 is the establishment of local weak comparison, which is done in \S \EE.


Note that the cases    where  $\psi$ takes values in the
{\bf interior} of the range of $f$,  
are much easier 
 than the general case  considered here.  Under this assumption the linearization of the operator is ofter
 quite nice.  In fact sometimes these cases
can be handled by results in [\DDR] (see Example 6.15 and results in \S 18).

Let us say again that in our past work  we have not considered operators, but rather we have
 taken a potential theory approach where the differential equation is given by 
 the boundary of a subequation.  Our main reason for considering 
 operators here is that with our hypothesis on $F$ and $f$ we can
 solve for general  inhomogeneous terms $\psi$.

An outline of our results is the following (details appear in the next section). 
We begin with a manifold $X$ and a pair $(F,f)$
where $f$ is an operator and $F$ is its domain.  That is, $F$ is a closed subset
of the 2-jet bundle   of $X$ and $f\in C(F)$.  We require $F$ to have the natural properties
of a  subequation as defined in [\DDR].  For the constant coefficient case in $\rn$
subequations are defined in Definition \AA.1.  For general equations considered here,
there are local automorphisms of the 2-jet bundle -- the jet-equivalences defined in 
Definition \AA.2 -- which take $(F, f)$ onto a constant coefficient  pair $(\bbf,  \bdf)$
in local coordinates.  The properties we need for $(\bbf,  \bdf)$ pull-back to
the desired properties for $(F, f)$.  In particular, $F$ is a subequation in the sense of [\DDR].
The local operator $\bdf$ is assumed {\sl tame} (Def. \AA.3) and compatible with $\bbf$
(Def. \AA.4).  We have tried to write this paper with a minimum of global geometry,
to reach a wider audience.  The global viewpoint is carefully presented in  [\DDR].
The main new part here is the local weak comparison   Theorem \EE.2 which
is a local result.

Now given such a subequation-operator pair $(F,f)$ and a function $\psi : X \to \bbr$ with values in the
range $\bdf(\bbf)$ of the operator, we want to solve the problem
$$
f(J^2 u ) \ = \ \psi \qquad{\rm with}\ \ J^2 u\in F\ 
\eqno{(\II.1)}
$$
at all points of a domain $\O\ss \ss X$ with prescribed continuous boundary values $\vf\in C(\bo)$.

At this level of generality, with no convexity or non-degeneracy assumption,
there is only one way available to give meaning to the equation (\II.1), 
namely, one of the equivalent viscosity definitions.
(See [\BEL] for the equivalence of the distributional approach when convexity
is assumed.)  To do this we consider the subset
$$
F_f(\psi) \ \equiv\ \{J\in F : f(J) \geq \psi\}
\eqno{(\II.2)}
$$
From our assumptions on $f$ (i.e., on $\bdf$) we see that the solutions
to our problem are  solutions to the $F_f(\psi)$-harmonic Dirichlet problem
as in [\DDR].  Utilizing Dirichlet duality, such a solution is a continuous function on $\ob$ such that
in $\O$
$$
\text{$u$ is $F_f(\psi)$-subharmonic \and $-u$ is $\wt{F_f(\psi)}$-subharmonic.}
$$
This means the following.
A continuous function $u$ on an open set $\O$ is $G$-{\sl subharmonic} for a subequation $G$ if for all
$x\in\O$  and all $C^2$-functions $\vf$ near $x$ with $u\leq \vf$
and $u(x)=\vf(x)$, we have $J_x^2\vf \in G$.  
If $G$ is a subequation, so is its {\bf dual}  $\wt G \equiv -(\sim \Int G) =\  \sim (-\Int G)$.
Under our assumptions here, $F_f(\psi)$ is a subequation, and one computes  that the dual is
$$
\wt{F_f(\psi)} \ =\ \wt F \cup \{J ; -J\in \Int F \ \ {\rm and}\ \ f(-J)\leq \psi\}.
\eqno{(\II.3)}
$$

The general pattern of our proof is to show (in \S \EE) that $F_f(\psi)$ satisfies
local weak comparison (Def. \EE.1). It then satisfies global weak comparison by [\DDR, Thm.\ 8.3].
Now of course some global hypothesis is required.
If there is a global approximator, then $F_f(\psi)$ satisfies comparison by
[\DDR, Thm.\ 9.7].
Theorem 5.2 in [\DDR] then gives the Main Theorem \AA.11.

To get a global approximator  we assume that $F_f(\psi)$ has a monotonicity cone $M$
(coming from one for the constant coefficient model).  We then assume that there exists
a smooth strictly $M$-subharmonic function on a neighborhood of $\ob$.

In  the general manifold case, such a function is certainly necessary.  Consider the inhomogeneous complex
Monge-Amp\`ere equation on a domain $\O$ in a complex manifold, and suppose that
it is always solvable as above.  Now blow-up a point  $x_0\in \O$ and choose a 
function $\psi$ which is positive on $D\equiv \pi^{-1}(x_0)$ where $\pi:\wt\O \to \O$
is the blow-up projection.  The Dirichlet problem is not solvable for this $\psi$.
This follows since any pluri-subharmonic function $u$ will be constant on $D$,
and hence the determinant of its complex hessian will be $\leq 0$ (actually $\equiv 0$ if $n>2$) along $D$.

Because there are so many important special cases of our Main
Theorem \AA.11 which are of historical significance in the literature, many examples and historical
remarks are given in Section \CC.  However we give a few examples just below to give 
an idea of the scope of the Main Theorem.

In Section \CCC\ we consider the case of solving the inhomogeneous equation with
 a measure $\mu$ on the right hand side.  This is sometimes possible with $\mu$
 taken to be the Dirac delta function.  However, in this case one needs the operator to be
 homogeneous and one must properly adjust its 
 homogeneity.
 
Now may be a good time for the reader to see the kind of  examples to which
our Main Theorem applies.  These and many more are treated in detail in Section \CC.

\Ex { \II.1. (The Monge-Amp\`ere operator on Almost Complex Manifolds)}
On any almost complex manifold $(X, J)$ there is an intrinsic operator $i \partial\dbar$
which allows one to define a subequation $F\equiv \cp_\bbc \ss J^2(X)$ consisting  at $x\in X$
 of $J^2_x u$ with $(i \partial\dbar  u)_x \geq 0$.  This allows us to define
the homogeneous complex  Monge-Amp\`ere equation by the boundary of $\cp_\bbc$, i.e.,
the $\cp_\bbc$-harmonics.

Now given a volume form $\O$ on $X$, we can define a Monge-Amp\`ere operator
$$
f(J^2 u) \ =\ { (i \partial\dbar  u)^n  \over \O} 
$$
($\dim_\bbr X = 2n$).  This gives an operator pair $(\cp_\bbc, f)$ and for any
continuous function $\psi\in C(X)$ with $\psi\geq 0$ we have the inhomogeneous equation
$$
f(J^2 u) \ =\ \psi. 
\eqno{(\II.4)}
$$
It follows from our Main Theroem (and it was already shown in [\AC])
 that the Dirichlet problem for (\II.4) can be solved for arbitrary continuous 
 boundary data on any compact, smooth domain with a strictly $J$-psh defining function.

\Ex { \II.2. (Invariant Operators on Riemannian Manifolds -  e.g., Krylov / Donaldson operators)}
On any riemannian manifold $X$ there is a  riemannian Hessian
operator on $C^2$-functions $u$ given for vector fields $V, W$ by
$$
(\Hess \, u)(V, W) \ =\ VW u - (\nabla_V W) u
\eqno{(\II.5)}
$$
 where $\nabla$ is the Levi-Civita connection on $X$.
 This Hessian is a symmetric
 tensor in $V$ and $W$, and gives a projection $J^2(X)\to \Sym (T^*X)$.
 
  Now  given an O$(n)$-invariant subequation $\bbf\ss\Symn$
 and an  O$(n)$-invariant operator $\bdf \in C(\bbf)$, then these give rise to a subequation
 $F$ and operator $f$ on $X$.  (This is well explained in [\DDR].)  To see that $F$ is a subequation
 it is only necessary to show that $\bbf+\cp\ss\bbf$ where $\cp\ss\Symn$ is the set of $A\geq 0$.
 The operator $f$ is tame and compatible (Def.'s \AA.3 and \AA.4) if $\bdf$ is.
 
 As an example consider the $k^{\rm th}$  Hessian operator given on $A\in\Symn$ by 
 ${\bf \s}_k(A) =   \s_k(\l_1, ... , \l_n)$, the $k^{\rm th}$ elementary symmetric symmetric function of
 the eigenvalues of $A$.  The natural domain for this operator is
 $$
 {\bf \Sigma}_k \ \equiv\ \{A\in\Symn : \s_1(A)\geq0, ... , \s_k(A)\geq 0\}.
 $$
 This pair has been studied for domains in $\rn$ by a number of authors (e.g., [\CNS], [\TRUU], [\TRU], [\TWC], 
 [\TWCC], [\TWCCC],  [\SPR]).
 Note that $k=1$ gives the riemannian Laplacian, and $k=n$ gives the riemannian real Monge-Amp\`ere 
 operator.
 
 Associated to these are the quotients
 $$
\s_{k,\ell} \ =\ {\s_k\over \s_\ell} \qquad{\rm on}\ \  {\bf \Sigma}_k
\eqno{(\II.6)}
 $$
 for $\ell < k$, which were studied by Krylov in [\KRY] and many others (see Spruck for example [\SPR]).
  Our Main Theorem \AA.11 solves the inhomogeneous Dirichlet Problem for this equation
 on manifolds.

 For $(k,\ell)=(n, n-1)$ these equations have received  much attention due to a conjecture of Donaldson
 (see [\DON]).

\Ex { \II.3. (Operators on G$_2$-manifolds)}
Let $X$ be a riemannian 7-manifold with G$_2$-holonomy (or more generally with a 
topological G$_2$-structure [\DDR, Ex.\ I in \S 1]).  Let $\bbf \ss\Sym(\bbr^7)$ be the set of 
$A$ with $\tr(A\bigr|_W)\geq 0$ for all associative 3-planes $W$.  Let $\bdf$ be the operator
$$
\bdf(A) \ \equiv\ \min\left\{ \tr\left(A\bigr|_W\right ) : W\ \text{an associative 3-plane} \right\}
$$
Then this gives a tame and compatible pair $(F,f)$ on $X$ to which the Main Theorem applies.

There is a similar story for the coassociative case.

\Ex { \II.4. (The Lagrangian Monge-Amp\`ere operator on Gromov Manifolds)}
Let $(X, \o)$ be a symplectic manifold equipped with a Gromov metric.  Set ${\bf Lag} \equiv
\{ A : \tr(A\bigr|_W)\geq 0$ for all Lagrangian planes $W\}$, and let ${\rm Lag} \ss J^2(X)$ be the subequation
determined as in \II.4.  The authors showed in [\LAG] that there is a natural polynomial differential
operator ${\bf M}_{\rm Lag}$ on ${\bf Lag}$, called the Lagrangian Monge-Amp\`ere operator.  It is tame and compatible 
with ${\bf Lag}$.  Thus this gives a natural operator ${\rm M}_{\rm Lag}$ on ${\rm Lag}$ to which our Main Theorem applies.
See Example \CC.7 and Theorem \CC.8 below.

Our Main Theorem \AA.11 has a generalization (Theorem \AA.11$'$) where 
the assumption of jet-equivalence is expanded to affine-jet-equivalence.

\Ex{\II.5}  This generalized Theorem \AA.11$'$ gives solutions to the Dirichlet
problem
$$
\det\left\{ \Hess_x u + M_x\right\} \ =\ \psi(x)
$$
on a riemannian manifold, where $M$ is a section of $\Sym ( T^*X)$.

Finally we recall the basic concept used for uniqueness in the Dirichlet problem. 
 Let $G$ be a subequation on a manifold $X$, and consider a domain $\O\ss\ss X$.
 By $G(\ob)$ we mean the set of upper semi-continuous functions on $\ob$ which are
 $G$-subharmonic on $\O$.

\noindent
{\bf Definition \II.6}
 We say that {\bf  comparison} holds for $G$ on $X$ if  for all $\O\ss\ss X$,
and for all $u\in G(\ob)$,  $v\in \wt G(\ob)$, one has that
$$
u+v \ \leq\ 0 \qquad {\rm on} \ \ \partial \O
\qquad\Rightarrow\qquad
u+v \ \leq\ 0 \qquad {\rm on} \ \  \ob.
$$
 
 Note that if $u$ and $w$ are solutions to the Dirichlet problem on $\O$,
 then $u,w\in  G(\ob)$, $-u, -w \in \wt G(\ob)$ and $u=w$ on $\bo$. Hence, comparison implies 
 that $u=w$.

\Note{\II.7}  Of course an interesting case of the work here is when
$(F,f)=(\bbf,\bdf)$ is itself    constant coefficient in euclidian space.
This case (pure second-order) is contained in the work of Cirant and Payne [\CP], where
other quite nice theorems are proved.

\vskip .3in

\vfill\eject

\noindent {\headfont \AA.    Statement of the Main Result.}     

\medskip
We  begin by considering the constant coefficient (or euclidean) case.  
Let 
$$
\bbj^2 \ =\ \bbr\oplus \rn \oplus \Symn
\eqno{(\AA.1)}
$$
be the space of  ``2-jets at 0"   
with classical coordinates  $(r, p, A)$. 

\Def{\AA.1}  By a {\bf constant coefficient subequation} on $\rn$ we mean a closed subset  $\bbf\ss \bbj^2$
such that

(P)  \quad  $\bbf+ (0, 0, P) \ \ss\ \bbf \qquad \forall\, P\geq0$ \ \ \  (Positivity or Weak Ellipticity),

(N) \quad $\bbf+ (-r, 0, 0) \ \ss\ \bbf \qquad \forall\, r\geq0$\quad\ (Negativity),

(T)  \quad $\bbf \ =\ \overline{\Int \bbf}$\qquad \qquad \qquad\qquad \qquad   (The Topological Condition)

\Def{\AA.2} Given such an $\bbf$,  we consider 
 a continuous function, or {\bf operator}
$$
\bdf \in C(\bbf).
$$
We call $(\bbf, \bdf)$ a {\bf constant coefficient subequation-operator pair}
(often shortened to ``operator pair'' when the meaning is obvious).
Note that the case $\bbf  =  \bbj^2$ is allowed here.

We introduce the following structural condition on the operator $\bdf$.

\Def{\AA.3.  (Tameness)}  The operator $\bdf\in C(\bbf)$ is said to be {\bf tame}  
on $\bbf$ if 
$$
\begin{aligned}
&\forall\, s, \l > 0\ \ \ \exists\, c(s, \l)\ >\ 0\ \ \ 
{\rm such\  that\ \ } \\
 \bdf(J+ (-r, 0, P)) \  - \ & \bdf(J)\ \geq\ c(s, \l) \quad \forall J \in \bbf,  r \geq s \ {\rm and}\  P\geq\l I.
 \end{aligned}
\eqno{(\AA.2)}
$$
This is a mild\footnote{(``mild'' in the sense that  it holds for most natural operators.
(See Propositions \CC.11 and \CC.13.)} strengthening of the required, weakest possible assumption:

\noindent
{\bf (Degenerate Elliptic) on $\bbf$}
$$ 
\bdf(J+(0, 0, P)) - \bdf(J) \ \geq\ 0 \qquad \forall\, J\in \bbf\ \ {\rm and}\ \ \forall \, P\geq0.
\eqno{(\AA.3)}
$$

Note that if $\bdf$  is  degenerate elliptic, then if condition (\AA.2) holds for $c=c(\l)$ and  $P=\l I$.
 It also holds for $P\geq \l I$
by (\AA.3).

\def\bbff{\bbf/\bdf}

There is a second condition we must impose, which is a compatibility between  
the operator $\bdf$ and the subequation $\bbf$
(when $\bbf$ is not all of $\bbj^2$).

\Def{\AA.4.  ($\bbff$-Compatibility)}  
We say that the set $\bbf$ and the operator $\bdf$ are {\bf compatible} if
$$
\partial \bbf \ =\ \{\bdf=c_0\}\qquad \text{for some }\ \ c_0\in\bbr.
\eqno{(\AA.4)}
$$
As explained in Proposition B.1 in Appendix B, this condition is necessary in our main theorem.
It implies that the level sets $\{\bdf=c\}$, for $c>c_0$ are contained in $\Int \bbf$,
i.e., they do not meet the boundary $\partial \bbf$. For instance, it
eliminates the following ``bad'' case.

\Ex {\AA.5} Consider the pure second-order  subequation on $\bbr^2$:
 $\bbf = \cpt = \{\l_{\rm max}\geq 0\}$, and let $\bdf=\l_1+\l_2$.  Here $\bdf(\bbf)=\bbr$, and
for all $c<0$ the boundary of $\bbf_c\equiv \{\bdf\geq c\}$ contains points of $\partial \bbf$ where $\bdf>c$.
There are lots of examples like this one, where $\bdf$ is elliptic on $\bbf$, but
$$
\text{
 $\exists\ \  c\in \bdf(\bbf)$ and $J\in\partial \bbf$ with $\bdf(J)>c$.}
\eqno{(\AA.4)^*}
$$
Note that (\AA.4)$^*$ is the negation of (\AA.4).

The final ingredient is the following.

\Def{\AA.6}  Let $(\bbf, \bdf)$ be a operator pair.  By a {\bf monotonicity cone} for $(\bbf, \bdf)$
we mean a constant coefficient convex cone subequation $\bbm\ss \bbj^2$, with vertex at 0,
such that 
$$
\bbf(c) + \bbm \ \ss\ \bbf(c) \qquad\text{for all values $c$ of $\bdf$}
$$
where $\bbf(c)\equiv \{\bdf\geq c\}$.

The pure second order case of the following result follows from the work of  Cirant and Payne [CP].
In fact their work is much more general; they consider operators of the form $f(x, D^2u)$.
For the cases considered here their assumptions are equivalent to our tameness condition.

\Theorem{\AA.7. (Constant Coefficient Operators)}
{\sl
Let $(\bbf,  \bdf)$ be a compatible  operator pair where $\bdf$ is tame on $\bbf$,
and suppose $\bbm$ is a monotonicy cone for $(\bbf,  \bdf)$.
Let $\O\ss\ss\rn$ be a domain with smooth boundary which 
satisfies the strict boundary convexity condition (Def.\ \AB.1. See also Thm. \AB.5.).
Suppose also that $\ob$ admits a smooth strictly $\bbm$-subharmonic function.
Then for each $\psi\in C(\ob)$ with values in $\bdf(\bbf)$,  and each  $\vf\in C(\bo)$,
there exists  a unique function $h\in C(\ob)$ satisfying:

\medskip
(1)\ \ $h$ is a   (viscosity) solution to $\bdf(J^2 u)=\psi, \ J^2 u \in\bbf$ on $\O$, and

\medskip
(2)\ \ $h\bigr|_{\bo} =\vf$.

\medskip\noindent
Furthermore, comparison holds, and  $h$ is the associated Perron function.
}

\medskip

Note that if $\psi \in C(\ob)$ does not take its values in $\bdf(\bbf)$, then problem (1.1) makes no sense
for smooth functions.
The functions $\psi \in C(\ob)$ which satisfy this necessary condition: $\psi(\ob) \ss \bdf(\bbf)$ will be 
called {\bf admissible (inhomogeneous terms)}.

\Remark {\AA.8}  The notion of strict $\bbf$   convexity for $\bo$ 
goes back to Caffarelli,   Nirenberg and Spruck
[\CNS], and appears in many works of the authors. 
The concept is  discussed   in Section \AB.

Suppose now that an operator $\bdf\in C(\bbf)$ has 
 the property that for some strictly increasing  continuous function
$\chi$ defined on the set $\bdf(\bbf)\ss \bbr$, the operator $\overline \bdf \equiv \chi\circ \bdf$ is tame on $\bbf$.
Then  $\bdf$ is said to be  {\bf tamable (by $\chi$)}.

\Theorem{\AA.7$'$}
{\sl The conclusions of Theorem \AA.4 remain true for any   operator $\bdf\in C(\bbf)$  which can be tamed.}

\pf    Set $\overline \psi = \chi\circ \psi$ and note that $\overline \psi$
is an admissible inhomogeneous term for $\overline \bdf \equiv \chi\circ \bdf$ if and only if $\psi$
is an admissible inhomogeneous term for $\bdf$.
\qed

\vskip .3in

\centerline
{
\bf  Second order equations on a manifold.
}

\medskip
We now take up the discussion of subequation-operator pairs $(F,f)$ on an $n$-manifold $X$.
We recall that the natural setting for second-order equations is the 2-jet bundle
$J^2 X \to X$ defined intrinsically at a point $x\in X$ as the quotient 
$$
J^2_x(X) \ \equiv\ C^\infty_x / C^\infty_{x,3}
$$
where  $ C^\infty_x$ denotes the germs of smooth functions at x, and $C^\infty_{x,3}$
the subspace of germs which vanish
to order three at x.  Given a smooth function $u$ on $X$, let  $J^2_x u \in J^2_x(X)$ denote its 2-jet at $x$, and
note that $J^2 u$ is a smooth section of the bundle $J^2(X)$.  This bundle is discussed in general in 
[\DDR].  However, we will only need the following.  Given a system of local coordinates  $\UU\ss\rn$ 
for  $X$, there is a natural trivialization
$$
J^2(\UU) \ =\ U\times (\bbr\oplus \rn\oplus \Symn)
\eqno{(\AA.5)}
$$
and $J^2_x u = (x, u(x), D_xu, D^2_xu)$.

The notion of jet-equivalence is crucial for this paper.  This concept is 
defined and broadly discussed on manifolds in [\DDR].  However, here we will only need  
to understand  it in the local trivialization (\AA.5).

\Def{\AA.9}  A (linear) {\bf jet-equivalence of $J^2(\UU)$} is a bundle automorphism
$$
\Phi : J^2(\UU) \ \arr\ J^2(\UU)
$$
given by
$$
\Phi(x,r,p,A) \ =\ (x, \ r, \ g(x)p, \ h(x) A h^t(x) +L_x(p))
$$
where $g,h: \UU\to {\rm GL}_n(\bbr)$ and $L:\UU \to \Hom(\rn, \Symn)$
are smooth (or at least Lipschitz continuous) functions.

\medskip

We point out that  jet-equivalences vastly change subequations.  For a smooth function $u$,
$\Phi(J^2u)$ is essentially never the 2-jet of a function.

\Def{\AA.10}  Let $F\ss J^2(X)$ be a closed set and $f\in C(F)$ an operator.  The the pair 
$(F,f)$ is {\bf locally jet-equivalent} to a constant coefficient operator pair $(\bbf, \bdf)$
if each point $x\in X$ has a local coordinate neighborhood $\UU \ss \rn$
and a jet-equivalence $\Phi : J^2(\UU) \to  J^2(\UU)$ which takes the pair
$(F,f)$ to $(\bbf, \bdf)$, that is,
$$
\Phi \left(F\bigr|_\UU\right) \ =\ \UU\times \bbf
\and
f \ =\ \bdf\circ \Phi.
$$
If, in addition, $\bbm$ is a monotonicity cone for $(\bbf, \bdf)$, and $M\ss J^2(X)$  is a closed set
such that for each local jet-equivalence above
$$
\Phi \left(M\bigr|_\UU\right) \ =\ \UU\times \bbm,
 $$
we say that $(F, f), M$ is {\bf locally jet-equivalent} to $(\bbf, \bdf), \bbm$.

\Theorem{\AA.11.\  (The Main Result)}
{\sl
Suppose that $(F, f)$ is a subequation operator pair with monotonicity cone  $M$
on a manifold $X$. Suppose further that
 $(F,f), M$  is locally jet-equivalent 
to a compatible  constant coefficient operator pair $(\bbf, \bdf)$ with monotonicity cone $\bbm$,
and that $\bdf$ is tame on $\bbf$. Let $\O\ss\ss X$ be a domain with smooth boundary, 
and assume the following given data:

 \noindent 
{\bf  Inhomogeneous Term:}
 
 (i)\ \ \  $\psi \in C(\ob)$ with values in ${\bf f}(\bbf)$, and

 \noindent 
{\bf  Boundary Values:}

 (ii) \ \ $\vf \in C(\bo)$.

If $X$ supports a smooth strictly $M$-subharmonic function, then comparison holds (Def. \II.6)
for arbitrary domains $\O\ss \ss X$.

If in addition  $\partial \O$ is smooth and  satisfies the strict boundary convexity condition (Def.\ \AB.1),
 there exists a unique $h\in C(\ob)$ which 
 
 (iii) \ \ satisfies the equation $f(J^2 h) = \psi$ on $\O$ (in the viscosity sense), and
 
 (iv) \ \ $h\bigr|_{\bo} = \vf.$

\medskip\noindent
Furthermore,   $h$ is the associated Perron function.
}

\Note{\AA.12}  
(a) For reduced subequations one can simply invoke strict $M$-convexity of $\bo$
 instead of using Def.  \AB.1 (see Theorem \AB.5 below).  
 
(b)  When the euclidean model $(\bbf,\bdf)$ is pure second-order,  the convexity cone 
subequation $P$,
with euclidean model  ${\bf P} = \bbr\oplus \rn\oplus \{A\geq 0\}$,
is always  a monotonicity cone for $(F,f)$ on $X$.  However for many such examples the optimal monotonicity
cone is much larger.

\medskip

 Theorem \AA.11 has a stronger version.

\Theorem{\AA.11$'$} {\sl
Theorem \AA.11 remains true if one replaces jet-equivalence with the more general concept of  affine jet equivalence
(see Def. \EE.4).
}

\vskip.3in


\noindent{\headfont \AB.\   Boundary Convexity}

The notion of boundary convexity of a domain is used classically to construct 
barriers, which are crucial in proving existence for the Dirichlet problem.
 Caffarelli,  Nirenberg and Spruck  [\CNS] presented a definition which worked for
 constant coefficient subequations in $\rn$, which are  orthogonally invariant and
 pure second-order.  Their ideas were adapted, first in [\DDD, \S 5] without any
 invariance, and then in [\DDR, \S 11] to the completely general case  of an
 arbitrary subequation on a manifold.  
 
The reader is referred to  \S 7 of [\SURVEY] for nice presentation of these ideas with many examples.

 \Def{\AB.1}  Let $\O\ss \ss X$ be a domain with smooth boundary in a manifold $X$.  
Let $(F,f)$ be an operator pair and $\psi\in C(\ob)$ an admissible inhomogeneous term as in Theorem 
 \AA.11. Then we have the subequation $F_f(\psi)$ and its dual given in (\II.2) and (\II.3).
We say that {\bf   $\bo$ satisfies the strict boundary convexity condition} if   each point $x\in\bo$
 is strictly 
 $F_f(\psi)$- and $\wt{F_f(\psi)}$-convex, as defined in \S 7 of [\SURVEY].
 
 Now in the case where $(F, f)$ is reduced (i.e., independent of the dependent variable), 
this condition is implied by a simple condition that depends only on the  monotonicity
 cone $M$.   To state this we recall some basic definitions and prove a  Lemma 
 which has some independent interest.
  
 \def\Jred{J^2_{\rm red}(X)}
 
 We recall that there is a canonical splitting $J^2(X) = \bbr \oplus  \Jred$
(where $\bbr$ corresponds to the value of the function).  By a reduced subequation
we mean one of the form $\bbr\oplus G\ss\bbr\oplus \Jred$.  For the rest of this section
all subequations will be reduced.

 Given a reduced subequation  $G\ss \Jred$  on a manifold $X$ there is an associated
 {\bf asymptotic interior} $\oa G$ where $J\in \oa G$ if there is an open set $J \in \cn(J)\ss \Jred$
and $t_0>0$ with 
$$
t\,  \cn(J) \ \ss\ G \qquad {\rm for\ all\ }\ t\geq t_0.
\eqno{(\AB.1)}
$$
 This defines an open set $\oa G\ss \Jred$ which is a bundle of cones with vertices
 at the origin in each fibre.
 
 It is immediate from this definition that for any two subequations
 $$
 G\ \ss\ H \qquad\Rightarrow \qquad \oa G\ \ss\ \oa H.
 \eqno{(\AB.2)}
$$

 Moreover, one see easily that (with vertices at the origin)
 $$
 \text{If $G$ is a cone subequation, then} \ \ \oa G \ =\ \Int G.
 \eqno{(\AB.3)}
$$
 
 The assertion can be carried over to translates as follows.
 
 \Lemma{\AB.2} {\sl
 Suppose that $G$ is a cone subequation with vertices at the origin, and $J_0$ is a continuous 
 section  of $\Jred$. Then the translated subequation has the same asymptotic interior:}
 $$
 \oa{{G+J_0}} \ =\ \Int G.
  \eqno{(\AB.4)}
$$
 \pf
 Suppose $J\in  \oa{{G+J_0}}$, i.e., there exists a neighborhood $\cn(J)$ and $t_0>0$ such that
 $t \cn(J) \ss J_0+ G$ for all $t\geq t_0$.    Then $ts\cn(J)\ss J_0+G$  for all $t\geq t_0$ and $s>1$.
 Since $G$ is a cone bundle, we have $t\cn(J) -{1\over s} J_0 \ss G$.  Sending $s\to \infty$ proves that
 $t \cn(J) \ss G$ for all $t\geq t_0$.  That is, $J \in \oa G$,  which by (\AB.3) equals $\Int G$.
 
 Conversely,  if $J\in \Int G$, then there exists $\cn(J) \ss \Int G$.  
 Since $\Int G$ is a bundle of cones, $t\cn(J)\ss\Int G$ for all $t>0$.
 Pick a small neighborhood  $\cn'(J)$ and $t_0>0$ so that  $\cn'(J) - {1\over t} J_0 \in \cn(J)$ for all $t\geq t_0$.
 Then   $t \cn'(J) \ss J_0 + t\cn(J) \ss J_0+G$ for all $t\geq t_0$ proving that $J\in \oa{{J_0+G}}$. \qed
 
 The interior of a monotonicity subequation for $G$ is smaller than the asymptotic interior of $G$.

 \Cor {\AB.3}  {\sl Suppose $M \ss \Jred$ is any bundle of cones with vertices at 0.  If $G\ss \Jred$ is a subequation which is $M$ monotone, then }
 $$
 \Int M \ \ss\ \oa G.
 \eqno{(\AB.5)}
$$ 
 \pf
 Fix $x\in X$ and choose a local section $J_0$ of  $\Jred$ defined near $x$ and taking values in $G$.
 Since $G$ is $M$-monotone, $J_0+M\ss G$.  Hence, by (\AB.2),  $\oa{{J_0+M}}\ss\oa G$. 
  Finally, by Lemma \AB.2, $\Int M = \oa{{J_0+M}}\ss\oa G$. \qed
  
  This extends as follows.
 
 \Cor{\AB.4}  {\sl
 Suppose $(F,f)$ is a reduced subequation pair and $M\ss \Jred$ is any  bundle  of cones  with vertices at 0. 
  If $F$ is $M$-monotone and the operator $f$  is $M$-monotone, then for each admissible $\psi$ the inhomogeneous subequation
 $F_f(\psi)$ is $M$-monotone, and hence
 }
 $$
 \Int M \ \ss\ \oa{{ F_f(\psi)}}
  \eqno{(\AB.6)}
$$
 
Suppose now that $(F, f)$ and $M$  is a  reduced triple, as above.  
 In this case the strict $F_f(\psi)$ convexity 
at $x\in\bo$, given in Definition \AB.1, is simply that in a neighborhood of $x$: \medskip

\vfill\eject

\centerline{ there exists a local smooth defining function
 for $\bo$}
 \centerline{ which is strictly $ \oa{{ F_f(\psi)}}$-subharmonic. }
 
 Now if the subequation $F_f(\psi)$ is $M$-monotone, so is its dual.  As a consequence we have the 
 following theorem.  We say that a boundary is {\bf strictly $M$-convex} if each point has a smooth local defining function
 which is strictly $M$-subharmonic, i.e.,  such that $J^2_{\rm red} \rho \in \Int M$.

 \Theorem{\AB.5} {\sl
 Let $(F,f)$ be a  operator pair with monotonicity cone $M$ as in Theorem \AA.11.
 If the triple $(F,f), M$ is  reduced, then any boundary which is strictly $M$-convex, satisfies the
 strict boundary convexity condition  \AB.1.
 }


 \vskip .3in


  \def\bbbf{{\bf f}}
 \def\Ffp{F_f(\psi)}
 \def\Fp{F(\psi)}

\noindent{\headfont \EE.\   Local Weak Comparison}

 Suppose that $G\ss J^2(X)$ is a subequation on a manifold $X$.
  Fix a metric on the 2-jet bundle $J^2(X)$.
 For $c>0$ we define $ G^c$ by its fibres
 $$
 G^c_x \ \equiv \ \{J\in G_x : \dist(J, \sim G_x) \geq c\}
 \ =\ \{J\in G_x  :  J+\eta \in G_x \ \ \forall\,\|\eta\|\leq c\}.
 $$

\noindent
{\bf Definition \EE.1}  We say that {\bf weak comparison} holds for $G$ on an open set $Y\ss X$ if there is a $c>0$ such that
for all $u\in G^c(Y)$,  $v\in \wt G(Y)$ and for all $\O\ss\ss Y$
$$
u+v \ \leq\ 0 \qquad {\rm on} \ \ \partial \O
\qquad\Rightarrow\qquad
u+v \ \leq\ 0 \qquad {\rm on} \ \  \ob
$$
 i.e., the {\bf Zero Maximum Principal} holds for $u+v$.
 We say that {\bf Local weak comparison holds for $G$ on $X$} if
 every point has a neighborhood $Y$ on which weak comparison holds.

 \Theorem{\EE.2.  (Local Weak Comparison)}
 {\sl
Let $(\bbf, \bbbf)$ be a constant coefficient subequation with operator on $\rn$
which is both compatible and tame.   Suppose that $(F, f)$ is a subequation with operator which
is jet equivalent to $(\bbf, \bbbf)$ on a open set $X\ss\rn$.   Then for any admissible continuous 
inhomogeneous term $\psi$, weak comparison hold for the associated
 inhomogeneous subequation $\Ffp  \equiv \{J \in F : f(J) \geq \psi\}$ on $X$
 }
 
\pf
Let $\Phi : J^2(X) \to J^2(X)$ be the jet bundle isomorphism taking  $(F,f)$ to $(\bbf, \bbbf)$, that is
$$
\Phi(F) \ =\ \bbf
\and
f \ =\ \bbbf \circ \Phi.
$$
In terms of the canonical trivialization of $J^2(X)$ we have for $x\in X$ that 
$$
(r', p', A') \ \equiv\ \Phi_x(r, p, A) \ \equiv \ (r, \ g(x) p, \ h(x) A h(x)^t + L_x(p))
\eqno{(\EE.1)}
$$
The associated inhomogeneous subequation $\Ffp \ss J^2(X)$ has fibre over $x\in X$
$$
\begin{aligned}
\Ffp_x \ &\equiv\ \{J\in F_x : f(x, J)\geq \psi(x)\}   \\
&=\  \{ J : J' \equiv \Phi_x(J) \in \bbf \ \ {\rm and}\ \ \bbbf(J') \geq \psi(x)\}
\end{aligned}
\eqno{(\EE.2)}
$$
The dual subequation $\wt{\Ffp}$ has fibre at $y\in X$ given by:
$$  
\wt{\Ffp}_y \ =\ \ft \cup \{J : -J \in \Int F \ \ {\rm and}\ \ f(y, -J) \leq \psi(y)\}.
\eqno{(\EE.3a)}
$$
Moreover,
$$  
J\in \wt{\Ffp}_y  \qquad \iff\qquad
J' \ \equiv\ \Phi_y(J) \in \wt{\bbf_{\bbbf}(\psi)}_y
\eqno{(\EE.3b)}
$$
and
$$  
 \wt{\bbf_{\bbbf}(\psi)}_y  \ =\ 
\wt \bbf \cup \{ J' : -J' \in \Int \bbf \ \ {\rm and}\ \ \bbbf(-J') \leq \psi(y)\}  
\eqno{(\EE.3c)}
$$

Failure of weak comparison for $\Ffp$ on $X$ 
means   there exists $\O \ss \ss X$, $c>0$, $u\in \Ffp^c(\ob)$ and 
$v \in \wt{\Ffp}(\ob)$ such that:
$$
u+v \ \leq \ 0  \ \  {\rm on} \ \ \partial \ob,\quad {\rm but}\ \ \sup_{\ob}(u+v) >0.
$$
(i.e., the Zero Maximum Principle fails  for $u+v$ on $\O$).
We use the Theorem on Sums of [CIL], in the form [\DDR, Thm. C.1].
It says that there exist a point $x_0\in \O$, a sequence of numbers $\e\searrow 0$
with associated points $z_\e = (x_\e, y_\e) \to (x_0, x_0)$, and
2-jets:
$$
\a_\e\ \equiv\ (r_\e, p_\e, A_\e)\ \in \ \Fp_{x_\e}^c
\and
\b_\e\ \equiv\  (s_\e, q_\e, B_\e)\ \in\ {\wt \Fp}_{y_\e}
 \eqno{(\EE.4)}
 $$  
 (for simplicity, here and below, we denote $\Ffp^c$ by $\Fp^c$, $\Ffp$ by $\Fp$, etc.)
 with the following properties.
$$
 r_\e \  =\ u(x_\e),\qquad 
 s_\e\ =\ v(y_\e),  \and r_\e+s_\e = M_\e \  \searrow \  M_0\ >\ 0
 \eqno{(\EE.5)}
 $$  
$$
 p_\e\ =\ \frac{x_\e-y_\e} {\e}\ =\ -q_\e \and  \frac{|x_\e-y_\e|^2} {\e} \ \arr\ 0
 \eqno{(\EE.6)}
 $$  
 $$  \left(
 \begin{matrix}  
 A_\e & 0 \cr 0 & B_\e
 \end{matrix}
   \right)
 \ \leq\  \frac{3}{\e} \left( 
  \begin{matrix}   
 I & -I \cr -I & I\cr
 \end{matrix}
  \right).
\eqno{(\EE.7)}
$$
We employ the notations
$$
\a_\e'\ \equiv\ (r_\e', p_\e', A_\e')  \ \equiv \ \Phi_{x_\e}(\a_\e)
\and
\b_\e'\ \equiv\  (s_\e', q_\e', B_\e')   \ \equiv \ \Phi_{y_\e}(\b_\e).
\eqno{(\EE.8)}
$$
By (\EE.1) this can be rewritten as
$$
r_\e' \ =\ r_\e,\quad p_\e' \ =\ g(x_\e) p_\e, \quad A_\e' \ =\ h(x_\e)  A_\e h(x_\e)^t + L_{x_\e}(p_\e)
\eqno{(\EE.9)}
$$
$$
s_\e' \ =\ s_\e,\quad q_\e' \ =\ g(y_\e) q_\e, \quad B_\e'\ =\  h(y_\e) B_\e h(y_\e)^t  +  L_{y_\e}(q_\e). 
\eqno{(\EE.10)}
$$

\Lemma{\EE.3}  {\sl
There exist $P_\e\geq0$ for $\e>0$ small, such that:}
$$
\lim_{\e\to 0}\left\{ \a_\e' + \b_\e' + (-M_\e, 0, P_\e)\right\} \ =\ 0.
$$

\pf
The first component  is $r_\e'-M_\e+s_\e' = r_\e-M_\e+s_\e$ which equals  zero by (\EE.5).
The second component   is 
$$
p_\e'+q_\e'  \ = \ g(x_\e){(x_\e-y_\e)\over \e}
-g(y_\e){(x_\e-y_\e)\over \e} \ =\ \biggl (g(x_\e)-g(y_\e)\biggr){(x_\e-y_\e)\over \e}
$$
 which converges to zero as $\e\to 0$ by (\EE.6).
It remains to find $P_\e\geq0$ so that the third component  
 $A_\e'+B_\e' + P_\e$, converges to zero.

 Multiplying both sides in (\EE.7) by
$$
\left( \begin{matrix} h(x_\e) & 0 \cr 0 & h(y_\e)  \end{matrix}  \right) \ \ {\rm on\ the\ left\ and\ \ }
\left( \begin{matrix} h(x_\e)^t & 0 \cr 0 & h(y_\e)^t\cr    \end{matrix} \right) \ \ {\rm on\ the\ right}
$$
gives 
$$
\left( \begin{matrix} h(x_\e) A_\e  h(x_\e)^t & 0 \cr 0 & h(y_\e) B_\e h(y_\e)^t\cr   \end{matrix}  \right)   
 \ \leq\ 
 {3\over\e} \left( \begin{matrix}  h(x_\e)h(x_\e)^t & -h(x_\e)h(y_\e)^t \cr
  -h(y_\e) h(x_\e)^t & h(y_\e) h(y_\e)^t  \end{matrix} \right).             
$$
Restricting these two quadratic forms to  diagonal elements $(x,x)$ then yields
$$
\begin{aligned}
h(x_\e) A_\e  h(x_\e)^t + h(y_\e) B_\e h(y_\e)^t  \ &\leq\ 
 {3\over\e}\left[  h(x_\e)( h(x_\e)^t -  h(y_\e)^t)  -     h(y_\e)( h(x_\e)^t -  h(y_\e)^t)   \right]  \cr
&=\ {3\over\e}  ( h(x_\e) -  h(y_\e)) ( h(x_\e)^t -  h(y_\e)^t)   \cr
&\leq \   {\l \over\e} \left| x_\e-y_\e\right|^2\cdot I \qquad{\rm for\ some\ \ } \l>0.
\end{aligned}
$$
Thus we can define  $P_\e\geq 0$  by:
$$
h(x_\e) A_\e  h(x_\e)^t + h(y_\e) B_\e h(y_\e)^t  + P_\e \ =\ 
 {\l \over\e} \left| x_\e-y_\e\right|^2\cdot I.
\eqno{(\EE.11)}
$$
It now follows from   the definitions in (\EE.9) and (\EE.10) that
$$
\begin{aligned}
 A_\e' + B_\e'   + P_\e \ &=\  {\l \over\e} \left| x_\e-y_\e\right|^2\cdot I 
+       L_{x_\e}(p_\e)   
+       L_{y_\e}(q_\e).    \cr
 \end{aligned}
\eqno{(\EE.12)}
$$
However,
$$
\begin{aligned}
 \left| L_{x_\e}(p_\e) + L_{y_\e}(q_\e) \right|
&=\    \left| ( L_{x_\e} -  L_{y_\e}) \left (  {x_\e-y_\e\over \e}  \right)  \right| \cr  
  & \leq \| L_{x_\e} -  L_{y_\e}\|  {|x_\e-y_\e|  \over \e}  \cr  
 =\ & \ \ O\left( {|x_\e-y_\e|^2  \over \e} \right)
\end{aligned}
$$
Using (\EE.6) this shows that 
$$
A_\e' + B_\e'   + P_\e \ \cong \ {|x_\e-y_\e|^2  \over \e}\ \to\ 0\quad{\rm as\ } \e\searrow 0. \qquad\mathqed
\eqno{(\EE.13)}
$$

We now examine the notion of $c$-strictness.  
Note that the definition of weak local equivalence is independent of the choice of metric on $J^2(X)$.

We set some notation.  If $\a \equiv (r, p, A)$ and $\eta$ are 2-jets at $x$, let 
$\a' \equiv (r', p', A')\equiv \Phi_x(\a)$ and $\eta' \equiv \Phi_x(\eta)$.  Since $\Phi_x:J_x^2(X)\to J_x^2(X)$ 
is a linear isomorphism, we can define a norm $\|\a\|$ on $J_x^2(X)$ to be the 
euclidean norm $|\a'|$ of $\a'=\Phi_x(\a)$.

By the definition of $c$-strictness, we have
$$
\a\in F(\psi)_x^c \qquad\iff\qquad
\a+\eta \in F(\psi)_x \quad\forall\, \|\eta\|\leq  c.
\eqno{(\EE.14)}
$$
By (\EE.2) we then have (with notation as above) that the first half of (\EE.4)
can be rewritten as
$$
\a_\e\in F(\psi)_{x_\e}^c \quad\iff\quad
\a_\e'+\eta' \in \bbf \ \ {\rm  and}\ \ \bbbf(\a_\e'+\eta') \geq \psi(x_\e)\quad\forall\, |\eta'|\leq  c.
\eqno{(\EE.4a)'}
$$
Using (\EE.3), the second half of (\EE.4) can be rewritten as
$$
(i) \ \ \b_\e' \in \wt \bbf \qquad {\rm or}\quad (ii) \ \ -\b_\e' \in \Int \bbf \ \ {\rm and} \ \ \bbbf(-\b_\e') \leq \psi(y_\e).
\eqno{(\EE.4b)'}
$$

We are now ready to complete the proof.  For $\e>0$ small enough, condition (i): $\b_\e' \in \wt \bbf$
is ruled out as follows.  If (i) holds, then by definition of the dual,  $$-\b_\e'  \notin  \Int \bbf.$$

Define 
$$
\a_\e'' \ \equiv \ \a_\e' + (-M_\e, 0, P_\e).
\eqno{(\EE.15)}
$$
By positivity (P) and negativity (N), for the subequation $\bbf(\psi)_{x_\e}^c$ and the fact that 
$\a_\e' \in \bbf(\psi)_{x_\e}^c$, it follows that:
$$
\a_\e'' \in \bbf(\psi)_{x_\e}^c.
\eqno{(\EE.4a)''}
$$

Now since $-\b_\e' \notin \Int \bbf$, we have that 
$$
0\ <\ c \ \leq\ {\rm dist}(\a_\e'', -\b_\e') \ =\ |\a_\e'', + \b_\e'|
$$
which, by Lemma \EE.2, has limit 0 as $\e\searrow 0$.  This shows that condition (i) is not
possible, and we are left with condition (ii).

Again, by the definition of $c$-strict, we can rewrite (\EE.4a)$''$ as
$$
\a_\e'' +\eta' \in \bbf     \quad{\rm and}\quad
\bbbf(\a_\e''+\eta')\geq \psi(x_\e) \quad \forall\, |\eta'|\leq c.
\eqno{(\EE.4a)''}
$$
Combining this with
$$
(ii) \quad    -\b_\e' \in \Int\bbf \qquad{\rm and}\qquad \bbbf(-\b_\e') \leq \psi(y_\e)
\eqno{(\EE.4b)'}
$$
yields
$$
\bbbf(-\b_\e') -  \bbbf(\a_\e'' +\eta' ) \ \leq \ \psi(y_\e) -\psi(x_\e) \qquad \forall\, |\eta'|\leq c.
\eqno{(\EE.16)}
$$

We shall now show that (\EE.16) violates tameness.  With $k, \l>0$ small and fixed,  define
$$
\eta_\e' \ \equiv \ - (\b_\e' + \a_\e'') - (-k, 0, \l I).
$$
Then $|\eta_\e' | \leq c$ for $\e>0$ sufficiently small by Lemma \EE.2, and so (\EE.16) holds.
However
$$
\a_\e'' + \eta_\e' + (-k, 0, \l I) \ =\ -\b_\e',
$$
so by the tameness of $\bbf$ the left hand side of (\EE.16) is bounded below by
the constant $c(k,\l)>0$, independent of $\e\to 0$.   Thus for $\e>0$ small, we have
$$0< c(k, \l) \leq  \psi(y_\e) -\psi(x_\e),$$  which is a contradiction since $y_\e-x_\e \to x_0-x_0=0$ as $\e\to0$.
\qed
\medskip

Theorem \EE.2 can be generalized by expanding the notion of equivalence.

\Def{\EE.4}  By an {\bf affine jet equivalence} we mean an automorphism
$\wt \Phi: J^2(X) \to J^2(X)$ of the form
$$
\wt \Phi \ =\ \Phi + J
$$
where $\Phi$ is a (linear) jet equivalence and $J$ is a section of the bundle $J^2(X)$.
 
 Suppose now that we have a subequation $F$ which is  affinely jet-equivalent to a
 constant coefficient equation $\bbf$ on a coordinate chart $U$.  Then 
it is shown in Lemma 6.14 in [\DDR] that if 
$$
J\in F_x \qquad\iff\qquad \Phi_x(J)+J_x \in \bbf,
$$
 then
 $$
 J\in \wt F_x \qquad\iff\qquad \Phi_x(J)-J_x \in  \wt \bbf
 $$

 We now go to the proof above where the hypothesis of jet equivalence
 is replaced by affine jet equivalence.  Then the display (\EE.8) must be replaced
 by
 $$
 \a_\e' = \Phi_{x_\e} + J_{x_e}
 \and
 \b_\e' = \Phi_{y_\e} - J_{y_e}.
 \eqno{(\EE.8)'}
 $$
 Since $J_{x_e}-J_{y_e} \to 0$ as $\e\to 0$, the proof goes through in this case.
 This give the following.
 
 \Theorem{\EE.5}  {\sl Theorem \EE.2 remains true if one assumes, more generally,  that $(F, f)$ is 
affinely jet equivalent to $(\bbf, \bbbf)$ (rather than just jet-equivalent to $(\bbf, \bbbf)$).
 }

\vfill\eject


\noindent{\headfont \BB .\   Proof of the Main Theorem}

We shall use the following.

\Theorem{\BB.1.\ (Thm.\ 9.7 in [\DDR])} {\sl
Suppose $F$ is a subequation on a manifold for which local weak comparison holds. 
Suppose there exists a $C^2$ strictly $M$ -subharmonic function
on $X$ where $M$ is a monotonicity cone for $F$. Then comparison holds
for $F$ on $X$.}
 
Now  on $X$  we are considering the subequation $F_f(\psi)$.  By Theorem  \EE.2
local weak comparison holds for this equation.  We have hypothesized that 
there is a strictly $M$ subharmonic function where $M$ is a monotonicity cone subequation
for $F_f(\psi)$. (See Definition \AA.6.) Hence comparison holds for $F_f(\psi)$ on $X$ by Theorem \BB.1 above.

The Main Theorem \AA.11 is now a consequence of  the following. 

\Theorem{\BB.2.\ (Thm.\ 13.3 in [\DDR])} {\sl Suppose comparison holds for a subequation $F$ on $X$.
Then for every domain $\O\ss\ss X$ with smooth boundary  which is strictly $F$- and
$\ft$-convex, both existence and uniqueness hold for the Dirichlet problem.}

\medskip\noindent
{\bf Proof of Theorem \AA.11.}  Use Theorem \BB.1 for uniqueness.  Use Theorem \BB.2 and Theorem
\AB.5 for existence. \qed

\medskip\noindent
{\bf Proof of Theorem \AA.11$'$.}  
This is the same, but one uses Theorem \EE.5 to get local weak comparison.

\vfill\eject


\def\ff{{\bf f}}

\noindent{\headfont \CC.\  Applications and Historical Remarks.}

The main result, Theorem \AA.11, applies to many equations of classical interest. 
We note, however,  that in these cases the operators $f$ are 
 almost always concave (so that the constraint sets are convex).
 By contrast, here  $F$ is an arbitrary subequation. 
  Furthermore,  in the literature
 the  inhomogeneous term $\psi$ is often required to satisfy a strict inequality $\psi>c$ where here
Theorem \AA.11 applies to any $\psi\geq c$ where $c$ is the minimum admissible value.

Now Theorem \AA.11 concerns subequation-operator pairs on manifolds
with the property that they  are locally jet-equivalent
to constant coefficient pairs $(\bbf, \ff)$.  As noted in \S \II \ such equations arise in a very
natural way -- for example, on almost complex manfiolds, on riemannian manifolds, on manifolds with
a topological reduction of structure group to $G\ss {\rm O}(n)$, etc.  
(see [\DDR], \SURVEY], [\AC]).
This certainly applies to 
manifolds with integrable reductions  (i.e., special holonomy) such as K\"ahler manifolds, hyperK\"ahler
manifolds, G$_2$  and Spin$_7$ manifolds, etc. 

Of course, Theorem \AA.11 does not address regularity,
and in fact,  without further assumptions no regularity beyond continuity is possible.\footnote{\  For an arbitrary
continuous  function $w\in C(\bbr)$, the function $u(x_1,x_2,x_3) = w(x_1)$ is $\bL_2$-harmonic on $\bbr^3$.}

Quite a few of the classical elliptic operators fall under a much more general rubric: homogeneous polynomials 
$\ff:\Symn \to \bbr$ which are {\bf G\aa rding hyperbolic} with respect to the identity $I$
(meaning $\ff(tI+A)$ has all real roots for each $A\in\Symn$).
Here one takes $\bbf = \overline\G$  where $\G$ is the G\aa rding cone, defined as the connected component
of $\{\ff>0\}$ which contains $I$, and one requires  $\ff$ to be weakly elliptic on $\bbf \equiv \overline \G$.
In all such cases there are many other branches of the 
equation (see  [\HYP]).  
Each branch has a natural operator which is covered by Theorem \AA.11.

We point out that  for such   G\aa rding hyperbolic operators $\ff$, {\bf the  G\aa rding cone $\bbf=\overline \G$ 
is a monotonicity cone}  for the pair $(\bbf, \ff)$ and also for all the other branches of the equation.

The operators given in Examples \CC.1 -- \CC.10 below are all G\aa rding hyperbolic 
with respect to the identity and tame. 
We point out that if  $\ff$ is  G\aa rding hyperbolic w.r.t.\ $I$, so are the  derivatives  ${d^k\over d t^k}\ff(tI + A)\bigr|_{t=0}$ 
(See [\HYP, Cor. 2.23]).  In Proposition \CC.11 we prove tameness 
is equivalent to being elliptic on $\bbf \equiv \overline\G$,
for all  G\aa rding polymonial operators, and this, in turn, is equivalent to 
$\cp \ss \bbf = \overline \G$.

A second different approach associates to any subequation $\bbf\ss\Symn$ a 
{\bf canonical operator} $\ff\in C(\Symn)$, which is  defined and tame on all of $\Symn$, with 
$\bbf = \{\ff\geq0\}$.  This completely general procedure is described  below. 
As an example, for $\bbf =\cp = \{A\geq 0\}$ (real Monge-Amp\`ere subequation)
the canonical operator is $\l_1(A)$.

We then exhibit operators which are topologically tame but not tamable,
also ones which are tamable but not tame.

At the end we discuss the asymptotic interiors for these many examples.

All the subequation-operator pairs $(\bbf,\ff)$ discussed in this section are compatible (Def. \AA.4).

\Ex{\CC.1. (Real Monge-Amp\`ere)}  The principal branch of this equation is:
$$
\det\left(  D^2 u\right )\ =\ \psi \quad \text{with $u$ convex and } \psi\in C(\ob), \  \psi\geq0.
\eqno{(\CC.1)}
$$
There is a long history of work on the principal branch beginning with the extensive work
of Alexandrov and Pogorelov. The reader is referred to Rauch-Taylor  [\RT] for a 
further discussion as well as a precise statement with two proofs.

Our main Theorem \AA.11 applies to the extension of this equation to any riemannian manifold $X$, namely
$$
\det\left(  \Hess \, u\right )\ =\ \psi \quad \text{with $u$ convex and } \psi\in C(\ob), \  \psi\geq0.
\eqno{(\CC.1)'}
$$
It asserts {\sl the existence and uniqueness of solutions to the Dirichlet Problem 
on any domain $\O\ss\ss X$ which supports a strictly riemannian convex function
and has a smooth strictly convex boundary} (the second fundamental form of $\bo$ with respect
to the interior normal is $>0$).

On the other hand, Theorem \AA.11 does not deal with the case where $\psi$ is a  measure,
which is done  in [\RT] when $X=\rn$.
Of course there are many results on this and related equations in $\rn$. See [\JTY] for a discussion and references.

The higher branches of this subequation are given in terms of the ordered eigenvalues $\l_1(A)\leq\cdots\leq \l_n(A)$ by
$$
\bL_k\ \equiv\ \{A : \l_k(A)\geq0\}.
\eqno{(\CC.2)}
$$
A tame operator for $\bL_k$, somewhat parallel to the determinant,  is given by
$$
\det_k(A) \equiv \l_k(A)\cdots\l_n(A).
\eqno{(\CC.3)}
$$
(However, for $k>1$ this is not a G\aa rding polynomial.)  The inhomogenous problem then becomes
$$
\det_k\left(  D^2 u\right )\ =\ \psi \quad \text{with $u\in \bL_k(\ob)$ and } \psi\in C(\ob), \  \psi\geq0.
\eqno{(\CC.4)}
$$

On the other hand the canonical operator associated to the $k^{\rm th}$ branch $\bL_k$
is just the $k^{\rm th}$ ordered eigenvalue function $\l_k$.
Since $\l_k$ is tame on all of $\Symn$, there is no restriction on the values
of the inhomogeneous term $\psi$.  Thus the inhomogeneous equation is given by
$$
\l_k\left(D^2 u\right) \ =\ \psi \qquad \text{with} \ \  \psi\in C(\ob).
\eqno{(\CC.5)}
$$
The Dirichlet problem for this equation was previously solved  in [\DDR] using the methods of local affine jet equivalence.

\Ex{\CC.2. (Complex  Monge-Amp\`ere)}  
The principal branch of this equation  in $\bbc^n$ is:
$$
\det_\bbc\left(  {\partial^2 u\over \partial z_i  \partial \bar z_j}  \right )\ =\ \psi \quad \text{where $u$ is psh and } \psi\in C(\ob), \  \psi\geq0.
\eqno{(\CC.6)}
$$
There is also a long history of work on this equation (usually under the assumption
that either $\psi=0$ or $\psi>0$).  The homogeneous case was 
initiated by Bremermann [\BRE] and then completed by Walsh in a short note [\WAL].
The solution in the inhomogeneous case was provided by the landmark paper 
of Bedford and Taylor [\BTA].  Since then many papers and books have added to 
this subject. 

This Dirichlet problem  was also solved on {\bf almost} complex manifolds
in [\AC] and [\PLI].  This is discussed in Example \II.1.

The higher branches  are treated exactly as in  (\CC.2) --- (\CC.4) 
except that one uses  the ordered eigenvalues of 
the hermitian symmetric matrix
$ \left(  {\partial^2 u \over \partial z_i  \partial \bar z_j}  \right )$.
Again one has the operator $\det_k$ as in (\CC.3). 
There is also the canonical operator $\l_k$,
degenerately elliptic on all of  $\Sym_{\bbr}(\bbc^n)$ as in  (\CC.5).

\Ex{\CC.3. (Quaternionic  Monge-Amp\`ere)}  The principal branch of this equation in $\bbh^n$  is:
$$
\det_\bbh\left( D^2 u  \right )_\bbh\ =\ \psi \quad \text{where  $u$ is \ $\bbh$-psh and } \psi\in C(\ob), \  \psi\geq0.
\eqno{(\CC.7)}
$$
By $A_\bbh$ we mean the quaternionic hermitian symmetric 
matrix ${1\over 4}(A-IAI-JAJ-KAK)$ whose eigenspaces are 
quaternion lines with eigenvalues $\l_1\leq\cdots \leq \l_n$, and $\det_\bbh A_\bbh\
\equiv \l_1\cdots \l_n$.
Results on the Dirichlet problem for this equation are due to Alesker [\ALL] and Alesker-Verbitsky [\AV].
However, there are higher branches of this equation,   defined in analogy with (\CC.2) --- (\CC.4),
to which our methods give new results.
Note that one  has two quaternionic operators, which are  analogues of (\CC.3) and (\CC.5).

\Ex{\CC.4. (The $k^{\rm th}$ Hessian Equation)}  

(a) {\bf The Real Case.} Consider the subquation ${\bf \Sigma}_k = \{A : \s_1(A)\geq0, ..., \s_k(A)\geq0\}$ where $\s_\ell$ denotes the $\ell^{\rm th}$ elementary symmetric function.
Now ${\bf \Sigma}_k$ is the closure of the connected component of $\{\s_k\neq0\}$
containing the identity $I$.  As with the previous examples there are $k$
(generalized) ordered eigenvalues, and therefore $k$ branches.
The principal branch of this  $k^{\rm th}$ hessian equation is 
$$
\s_k\left( D^2u  \right) \ =\ \psi  \quad \text{where  \ $u$ is \ ${\bf \Sigma}_k$-subharm. and } \psi\in C(\ob), \  \psi\geq0.
\eqno{(\CC.8)}
$$
This branch has been studied extensively by  Trudinger  [\TRUU] , [\TRU] and Trudinger-Wang 
 [\TWC], [\TWCC], [\TWCCC].

Of course using the riemannian hessian and our Main Theorem \AA.11, we have results on the Dirichlet problem 
for this equation on manifolds.

(b)  {\bf The Complex and Quaternionic  Cases.} Consider the analogous subquation ${\bf \Sigma}^\bbc_k = \{A : \s_1(A_\bbc)\geq0, ..., \s_k(A_\bbc)\geq0\}$ in $\bbc^n$ where $A_\bbc = \half(A-JAJ)$
is the hermitian symmetric part of $A$.  In analogy with (\CC.8) we obtain the  principal branch of the $k^{\rm th}$ complex hessian equation:
$$
\s_k \left(  {\partial^2 u\over \partial z_i  \partial \bar z_j}  \right ) \ =\ \psi  \quad \text{where  \ $u$ is \  ${\bf \Sigma}_k^\bbc$-subharm. and } \psi\in C(\ob), \  \psi\geq0.
$$
Work on this equation goes back to Blocki [\BL]. As in all other cases there are
branches and additional operators.

The quaternionic Hessian equation is the complete analogue of the example above with
 $A_\bbc$ replaced by $A_\bbh$.
 
 Theorem \AA.11 applies to solve the inhomogeneous Dirichlet problem for  these equations and their branches on manifolds.

\Ex{\CC.5. (The Quotient Hessian Equations)}
These are the operators $\s_{k,\ell} = \s_k/\s_\ell$ on ${\bf \Sigma}_k$
discussed at the end of Example \II.2.
Theorem \AA.11 applies to these are their complex and quaternionic analogues.

We finish the list of G\aa rding operators with several cases, which are  non-classical 
operators even for the principal branch.

\Ex{\CC.6. (The $p^{\rm th}$ Plurisubharmonic Equations)}  

(a)  {\bf The Real Case in $\rn$.} Consider the subequation 
$\cp_p \equiv\{A : \l_1(A)+\cdots+ \l_p(A) \geq0\}$ where  $\l_1(A)\leq\cdots\leq \l_n(A)$
are the  ordered eigenvalues of $A$.  Here there is a natural polynomial operator
$$
\det(\L^pA) \ \equiv\ \prod_{i_1<\cdots<\l_p}\left( \l_{i_1}+\cdots+\l_{i_p}\right)
$$
and an associated inhomogeneous equation for the principal branch:
$$
\det(\L^p D^2u)\ =\ \psi  \quad \text{where  \ $u$ is \ $\cp_p$-subharm. and } \psi\in C(\ob), \  \psi\geq0.
\eqno{(\CC.9)}
$$
This homogeneous Dirichlet problem for this equation  was solved in [\GEO, Theorem 7.6].
There are also  ${n\choose p}$ branches with operators defined in exact analogy
with the construction in Example \CC.1. This is obtained by using the (generalized)  eigenvalues $\l_I \equiv \l_{i_1}+\cdots+\l_{i_p}$. 

 Theorem \AA.11 applies to the inhomogeneous Dirichlet problem for this equation and its branches on manifolds,
 where the operator for the $k^{\rm th}$ branch
is the $k^{\rm th}$ ordered eigenvalue $\l_I$.

 \medskip
 
 (b)  {\bf The Complex Case in $\bbc^n$.} This is left to the reader.  It  parallels  the real case
 using the eigenvalues $\l_1,..., \l_n$ of $A_\bbc = \half(A-JAJ) \cong 
 \left(  {\partial^2 u\over \partial z_i  \partial \bar z_j}  \right ) $.

 \medskip
 
 (c)  {\bf The Quaternionic Case in $\bbh^n$.} This also  parallels  the real case,
but starting  with  the eigenvalues $\l_1,..., \l_n$ of $A_\bbh$.

\Ex{\CC.7. (The Lagrangian Plurisubharmonic Equation)}  Consider the subequation ${\bf Lag}$ in $\bbc^n$ defined by  requiring that $\tr\left(A\bigr|_W \right)\geq0$ for all Lagrangian $n$-planes $W$.
There is a U$(n)$-invariant polynomial operator ${\bf M}_{\rm Lag}$ defined on ${\bf Lag}$.
It depends only on the trace and the skew-hermitian part of the hessian, and it is a Lagrangian 
counterpart of the complex Monge-Amp\`ere operator.  This subequation and operator carry over
to any symplectic manifold equipped with a Gromov metric. All this is discussed in detail in [\LAG].
From Theorem \AA.11 we obtain the following.

\Theorem{\CC.8} {\sl
Let $X$ be a symplectic manifold with a Gromov compatible metric.  Suppose $\O\ss\ss X$ is a domain
with strictly Lag-convex boundary, which supports a strictly Lag-plurisubharmonic function.
Then for every continuous $\psi\geq 0$ on $\ob$ and every $\vf\in C(\bo)$, there is a unique 
 Lag-plurisubharmonic function $u$, continuous on $\ob$,  with
 
    (a)\ \ \   ${\rm M}_{\rm Lag}(u) \ =\ \psi$ \qquad\text{in the viscosity sense, and}

    (b) \ \ \  $u\bigr|_{\bo} \ =\ \vf$.
    }
 
 There are also results for the branches of ${\rm M}_{\rm Lag}$.

In all the following  examples we discuss the euclidean models.  However, the subequations and operators transfer 
to manifolds as discussed in Example \II.2, and Theorem \AA.11 applies.

\Ex{\CC.9. (The $\d$-Uniformly Elliptic Equation)}   The G\aa rding operator
$$
\ff_\d(A) \ \equiv\ \prod_{j=1}^n\left( \l_j(A) + \d \tr A    \right)
$$
on the principal branch (the G\aa rding cone) $\bbf\equiv \cp(\d) \equiv
\{\l_{\rm min}(A) + \d \tr A\geq0\}$
determines a uniformly elliptic inhomogeneous equation
$$
\ff_\d(D^2 u) \ =\ \psi \quad {\rm where} \ \ u\in \bbf(\O)\ {\rm and} \   \psi\in C(\ob), \ \ \psi\geq 0,
\eqno{(\CC.10)}
$$
to which Theorem \AA.11 applies.  All of the corresponding branches are also
uniformly elliptic, and Theorem \AA.11 applies similarly to them.
Of course Theorem \AA.11 also applies to their transfer to riemannian manifolds.

See theorem 5.16 in [\HYP] for a generalization with the eigenvalues 
$\l_j(A)$ replaced by the G\aa rding eigenvalues $\l^\ff_l(A)$ of 
an {\sl elliptic} G\aa rding operator as defined below.

\Ex{\CC.10. (The Pucci/G\aa rding Equation)}  This  is another G\aa rding operator related to the
standard {\bf  Pucci extremal operator} $\cp^-_{\l, \L}$,  which is defined for fixed constants $0<\l<\L$ by
$$
\cp^-_{\l,\L} \ \equiv\   \l \tr (A^+) +\L \tr(A^-) 
$$
where $A=A^++A^-$ is the composition of $A$ into $A^+ >0$ and $A^-<0$.
Associated to this is the subequation 
$$
{\bf P}_{\l,\L} \equiv \{ \cp^-_{\l,\L} \geq 0\},
$$
for which $\cp^-_{\l,\L}$ is the canonical operator (see Prop. \CC.13).
The monotonicity condition $\bbf + {\bf P}_{\l,\L}\ss\bbf$ is one of the 
many equivalent conditions of uniform ellipticity for a subequation $\bbf$.
Another is $\bbf + \cp(\d)\ss\bbf$. 
(See [\SURVEY, \S 4.5] for a detailed discussion.) 

Now we can define the  Pucci/G\aa rding polynomial $\ff_{\l,\L}:\Symn\to\bbr$ for which 
${\bf P}_{\l,\L}$ is the closed G\aa rding cone. It is constructed  using
the polar cone to ${\bf P}_{\l,\L}$, which  is the cone on $B_{\l.\L}\equiv \{\l I \leq A\leq \L I\}$.
The polynomial $\ff_{\l,\L}$ is then the product of the linear functions corresponding
to the vertices of the ``cube'' $B_{\l, \L}$.  
  This G\aa rding polynomial, which is of degree $2^n$,  can be explicitly computed.  
The minimum G\aa rding eigenvalue of $A\in\Symn$ is 
$ \cp_{\l,\L}^-(A) \equiv \l\tr A^+ + \L\tr A^- $.
Now $\cp_{\l,\L}^-(A)$ is customarily referred to as one of  the two  {\bf Pucci extremal operators} --
  the other being $\cp^+_{\l,\L}(A) = \l \tr A^- + \L \tr A^+$ which yields the largest  G\aa rding eigenvalue
   $\l\tr A^- +\L\tr A^+$.  Note that  the degree
of $\ff_{\l,\L}$ is high compared to that of $\ff_\d$, which is $n$. We refer to the
polynomial operator $\ff_{\l,\L}$ as the {\bf G\aa rding-Pucci operator} and the  equation
$$
\ff_{\l,\L}(D^2u) \ =\ \psi \quad {\rm where}\ u\in {\bf P}_{\l,\L}(\O) \ \ {\rm and}\ \psi\in C(\ob), \ \psi\geq0.
\eqno{(\CC.11)}
$$
as the {\bf inhomogeneous G\aa rding Pucci-equation}.
This equation and its branches make sense on any riemannian manifold, and Theorem \AA.11 applies.

We now make some general remarks.

\vskip .3in

\centerline{\bf Elliptic G\aa rding Operators}

Suppose $\ff$  is a   G\aa rding polynomial  on $\Symn$ of degree $m$,
which is $I$-hyperbolic.
The closed  G\aa rding cone $\bbf \equiv \overline \G$ is a convex cone  and as such is a subequation
if and only if $\cp\ss\bbf$.
In this case the operator $\ff$ is elliptic on $\bbf$.
This follows from the general fact that the G\aa rding eigenvalues $\l_j(A)$ are monotone precisely 
in $\bbf \equiv \overline \G$ directions.  Thus $\overline \Gamma$ is a monotonicity cone for 
each branch  $\bL_k = \{\l_k(A) \geq 0\}$ where $\l_1(A) \leq \l_2(A)\leq \cdots$ are the ordered 
eigenvalues.  (The reader is referred to   [\HYP] for a detailed discussion.)

\Prop{\CC.11} {\sl Each   G\aa rding polynomial $\ff$  with (closed)  
G\aa rding cone $\bbf = \overline \G \supset \cp$ is a tame operator
on $\bbf = \overline \G$  (its  principal branch). This pair $(\bbf, \bdf)$ determines
a pair $(F,f)$ on any riemannian manifold, and Theorem \AA.11 applies 
to the  inhomogeneous equation
$$
f(D^2 u)=\psi \quad {\rm where}\ u\in F(\O) \ \ {\rm and}\ \psi\in C(\ob), \psi\geq0.
\eqno{(\CC.12)}
$$
More generally, the operator $\ff_k(A) \equiv \l_k(A)\cdots \l_m(A)$ on the $k^{\rm th}$
branch $\bL_k$ is also tame, and has monotonicity cone $\bbf = \overline \G$.
Therefore Theorem \AA.11 applies to the extension of $(\bbf_k, \ff_k)$ to riemannian manifolds.
}

\pf
We must verify (\AA.2). Note that the ordered $\ff$-eigenvalues satisfy
$$
\l_k(A+\l I) \ =\ \l_k(A) +\l   \ \   {\rm if} \ A\in \Symn
\and
\l_k(A) \geq\ 0 \ \ {\rm if} \ A\in\bbf_k.
$$
Hence,
$\ff_k(A+\l I)-\ff_k(A) = \prod_{j=k}^m (\l_j(A)+\l) - \prod_{j=k}^m \l_j(A)  \geq \l^{m-k}$.\qed

\Remark{\CC.12}  It is easy to see that in all of the previous examples
one has $\cp \ss \bbf = \overline\G$ ( or equivalently $A\geq 0 \Rightarrow \l_k(A)\geq0$)
so that the  G\aa rding polynomial  $\ff$  is degenerately elliptic  on $\bbf$.
Consequently, by Proposition \CC.11 our main result  Theorem \AA.11 covers all of the operators in the first
ten examples above.

\vskip .3in


\centerline{\bf  Canonical Operators}

There is a canonical procedure for constructing an operator $\ff$ for an arbitrary subequation
$\bbf$.

\Prop{\CC.13} {\sl
For each subequation $\bbf\ss\Symn$ with $\bbf \ne \emptyset, \Symn$, and each normalizing
constant $k>0$, there exists a unique operator $\ff\in C(\Symn)$ satisfying
$$
\ff(A+\l I) \ =\ \ff(A) + k\l \qquad \forall\, A\in\Symn \ \ {\rm and}\ \ \forall\,\l\in \bbr
\eqno{(\CC.12)}
$$
and such that 
$$
\bbf \ =\ \{\ff(A)\ \geq\ 0\} \and \partial \bbf \ =\ \{\ff(A)\ =\ 0\}.
\eqno{(\CC.13)}
$$
Moreover, $\ff$ is tame so that Theorem \AA.11 applies to the inhomogeneous
Dirichlet problem}
$$
\ff(D^2 u) \ =\ \psi \qquad \psi \in C(\ob).
\eqno{(\CC.14)}
$$
\pf
The operator is constructed as follows. Consider the orthogonal splitting
$\Symn = \{\tr A=0\} \oplus \bbr\cdot I$ and choose coordinates $(x,y)$
($x=A-{1\over n}(\tr A)I, y={1\over n}\tr A$) with respect to this splitting.
Then there is a unique function $g(x)$  with the property that
$\bbf = \{(x,y) : y\geq g(x)\}$ and $\partial \bbf$ is the graph of $g$ over $\{\tr A=0\}$.
 The canonical operator $\ff$ is then defined by
$$
\ff(A) \ =\ ky - g(x)\ =\ \smfrac k n  \tr A     - g\left(A- \smfrac 1 n (\tr A) I  \right)
$$

This function $g$ is 1-Lipschitz with respect to norms
 $\|\cdot\|^\pm$ on $\{\tr A=0\}$ where $\|A\|^+ = -\l_{\rm min}(A)$ and 
$\|A\|^- = \l_{\rm max}(A)$.
See [\HYP, \S 3] (in particular Examples 3.4 and  3.5) for  details.
The proof that $\ff$ is tame is straightforward, with $c(\l) = k\l$.

\Remark{\CC.14} The two distinct methods of obtaining operators: (1) using a G\aa rding polynomial, and (2) constructing the canonical operator for a subequation can
be combined.  More precisely, given a subset $E\ss \bbr^m$ which is invariant
under permutation of coordinates and satisfies $E+\bbr^m_+\ss E$ (a ``universal eigenvalue subequation'') each degenerately elliptic G\aa rding operator $\ff$ of degree $m$ on $\Symn$ determines a new subequation $\bbf^\ff_E$ on $\rn$ by requiring that the G\aa rding eigenvalues of 
$A\in\Symn$ lie in $E$.  See Theorem 5.19 in [\HYP]  for details. If in addition $g$ is tame on $E$, adopting a straightforward definition, then $g(\l^\ff(A))$ is tame on $\rn$, and hence
Theorem \AA.11 applies to the inhomogeneous equation
$$
g\left(  \l^\ff(D^2 u)  \right) \ =\ \psi \quad {\rm where} \ u \in \bbf^\ff_E (\O)
\ {\rm and} \   \psi\in C(\ob), \ \ \psi\geq c
$$
where $c= \inf_E g$, and to the extension of this equation to riemannian manifolds.

\vskip .3in

\centerline{\bf Topological Tameness.}

\Def{\CC.15}  A function  $\bdf \in C(\bbf)$ is said to be {\bf topologically tame on $\bbf$} if
$$
\bdf (A+P) - \bdf (A) \ >\ 0 \qquad \forall\, A\in\bbf \ \ {\rm and}\ \ \forall\, P>0,
$$
or equivalently, if $\bdf $ satisfies ellipticity $\bdf (A+P) - \bdf (A) \geq 0$  and the above holds with $P\equiv \l I$, $\forall\, \l>0$.
The equivalence follows  since $P\geq \l I$ implies $\bdf (A+P)\geq \bdf (A+\l I)$.

 \Lemma{\CC.16}  {\sl Suppose that $\bdf $ is an elliptic operator on $\bbf$.  Then the following are equivalent:
 
 (1) \ \ The level sets $\{\bdf =c\}$ have no interior.
 
 (2) \ \ $\bdf $ is topologically tame on $\bbf$.
 
 \noindent
 In particular, all real analytic elliptic operators are topologically tame.}

\pf
If (2) is false, then for some $A\in\bbf$ and $P>0$ we have $\bdf (A+P)=\bdf (A)$ (using ellipticity).
Then for all $0<B<P$, we have $A+B\in \bbf$ and $\bdf (A+B)=\bdf (A)$ (by ellipticity).
This proves that $\{\bdf =c\}$ has interior where $c=\bdf (A)$, so (1) is false.

If (1) is false, pick $A\in \Int\{\bdf =c\}$.   Then $A+P \in \{\bdf =c\}$ for all $P>0$
sufficiently small proving (2) is false. \qed

\Cor{\CC.17} {\sl  Suppose $\bdf $ is an elliptic operator on $\bbf$ and comparison holds for the inhomogeneous
equation  when $\psi\equiv c$ is an admissible constant.  Then $\bdf $ must be topologically
tame.}

\pf
Suppose that the  level set $\{\bdf =c\} \equiv \{A\in \bbf : \bdf (A)=c\}$ has non-empty interior.  
Take  $u(x) = \half\bra {Ax}x$ so that $D^2 u= A \in\Int \{\bdf =c\}$.
 For all $C^2$ functions $v$ with compact
  support in a domain $\O$, and all $\e>0$ sufficiently small, 
  $$
  \bdf (D^2(u+\e v)) \ =\ c \qquad{\rm in}\ \ \O
  $$
  and $u$ and $u+\e v$ agree near $\bo$.  These counterexamples are all eliminated if
 the level set $\{\bdf =c\}$ has no interior.\qed

\vskip .3in

\centerline{\bf Non-Tame Operators.}

\Ex{\CC.18. (The Special Lagrangian Potential Equation / Topologically Tame but not Tamable for Certain Phases)} 
The operator $\ff\in C(\Symn)$ is defined by
$$
\ff(A) \ \equiv \ \tr \{\arctan(A)\}.
\eqno{(\CC.15)}
$$
Note that $\ff(\Symn) = (-{n\pi\over 2}, {n\pi\over 2})$.
This equation was introduced in [\CG] where it was shown that classical
solutions to  $\ff(D^2 u) = \theta$ have the property that the graph
of $Du$ in $\bbr^{2n}$ is special Lagrangian with phase $\theta$.
The important first work on this equation is due to Caffarelli, Nirenberg and Spruck [\CNS]
who established smooth solutions for $\theta$ in the outer-most branch where the 
subequation is convex. 
In [\DDD] existence and uniqueness were established for the continuous (DP) for $\ff(D^2 u) = \theta$
for all phases $\theta\in (-{n\pi \over 2}, {n\pi \over 2})$.  There is now
a copious literature.  Our purpose here is to discuss the inhomogeneous equation with $\psi(x)$ non-constant.
 For more historical comments the reader is referred to [\SLE].


\def\FTH{\bbf_{\Theta}}

\Prop {\CC.19}  {\sl
The degenerate elliptic operator $\ff(A) \ \equiv \ \tr \{\arctan(A)\}$ is topologically tame, but

(a) \ \  $
\bdf \  \text{ is {\bf not tamable} on $\FTH \equiv \{A : \bdf(A)\geq \Theta \}$ for $\Theta\leq (n-2){\pi\over 2}$.}
$

\noindent
However, for any

(b) \ \   $\Theta > (n-2)\pi/2$, the operator $\bdf$ {\bf is tamable} on the subequation $\FTH$}.

{\bf Proof.}  
Since $\bdf$ is real analytic, it is topologically tame.
Now  consider $A$ with 
$\l_1(A) <<0$ and $\l_k(A)>>0$ for $k>1$.  We can always choose these values
so that $\bdf(A) = (n-2){\pi\over 2}$.  As the absolute value of the eigenvalues becomes
very large the derivative of $\bdf(A)$ goes to zero.  Hence, no matter which smooth function $\chi$
one chooses, the composition $\chi\circ \bdf$ will have 
derivatives going to zero at these points, since $\chi'(\bdf(A))$ will not go to $\infty$
unless $\bdf(A)$ goes to ${n\pi\over 2}$. 

The proof of Part (b)  is given in [\SLE].  It was inspired by the result of Collins, Picard and Wu
 [\CPW]   that the subequation  $\FTH$ is convex for $\Theta > (n-2)\pi/2$,
 even though $\bdf$ is not concave unless $\Theta \geq (n-1){\pi\over 2}$.
\qed

A corollary of Proposition \CC.19 Part (b) is that comparison holds for the inhomogenous Dirichlet Problem
 $\bdf\left( D^2_x u\right) = \psi(x)$ on a domain $\O\ss\ss\rn$ provided that 
$\psi\in C(\ob)$ has values  $\psi(\ob) \ss  \left ((n-2){\pi\over 2}, n {\pi\over 2}\right)$
since our main Theorem  \AA.11' applies.
This comparison result was recently proved by  S. Dinew, H.-S. Do and T. D. T\^o  in [\DTT].

Existence for this Dirichlet problem requires computing the asymptotic cone for  the subequation $\FTH$.
For $\Theta > (n-1){\pi\over 2}$ this was done in [\CNS].  The main point of the article [\SLE]
is to compute this asymptotic cone for all $\Theta$, thereby providing the widest class of domains
$\O$ where existence holds.


Comparison  for a general admissible $\psi$, 
remains a difficult open question.  (See Question A in \S 8.)

\Ex{\CC.20. (Tamable but not Tame)}  Perhaps the simplest example is to start 
with the Laplace subequation $\bbf=\D\equiv \{\tr A\geq0\}$. Then the operator
$\ff(D^2u) \equiv \log(1+\tr D^2u)$ is not tame, since
$$
\ff(A+\l I) - \ff(A) \ =\ \log\left(  1+ {n\l \over 1+\tr A}\right)
$$
has infimum zero over $\tr A \geq c$.  However, $\chi(t) \equiv e^t -1$ tames $\ff$ since
$\chi\circ \ff(A) = \tr A$.

\Ex{\CC.21. (Another Topologically Tame but Non-Tamable Operator)} Define a topologically tame operator $\ff$ as follows.
First make the change of coordinates $$y=\tr A \and x=A-{1\over n}(\tr A) I.$$ 
Then set $\bbf  \equiv \{y\geq0\} = \{\tr A\geq0\} =\D$ 
and define a function $\ff\in C(\D)$ by
$$
\ff(x,y) \ =\ 
\begin{cases}
{y\over 1+\|x\|} \qquad\ \  {\rm if}\ \ y\leq 1+\|x\| \\
y-\|x\|    \qquad  {\rm if}\ \ y -\|x\| \geq 1. \\
\end{cases}
$$
{\bf Claim:} {\sl This operator $\ff$ cannot be tamed.}

\def\of{\bar \ff}

\pf Suppose $\of = \chi\circ \ff$ satisfies the tameness condition (\AA.6).
Choose $y>0$ and $\l>0$ small.  Then since $\of$ is constant on the level sets
of $\ff$   
$$
\begin{aligned}
\of(0,y+\l)\ &=\ \of(x, (1+\|x\|)(y+\l)) \and \\ \of(0,y)\ &=\ \of(x, (1+\|x\|)y)
\end{aligned}
\eqno{(\CC.16)}
$$
for all $x$.  Let $\|x\| = k \in\bbz^+$.  Applying condition (\AA.6) repeatedly shows that
$$
 \of(x, (1+\|x\|)(y+\l)) - \of(x, (1+\|x\|)y) \ \geq\ (1+\|x\|)c(\l) \ = \ (1+k)c(\l) \ \to\ \infty
$$
as $k\to\infty$.  However, by (\CC.16) we have
$$
\of(x, (1+\|x\|)(y+\l)) - \of(x, (1+\|x\|)y) \ =\ \of(0,y+\l) - \of(0,y).\qquad \mathqed
$$

\Ex{\CC.22. (Another Non-tamable Operator)} A similar, even wilder operator $\ff$
can be constructed on $\Sym(\bbr^n)$ as follows. We define $\ff$ in terms of the eigenvalues
of $A$ with the property that 
$$
\ff(A) \ =\ 
\begin{cases}
\l_{\rm min}(A) \qquad\ \  {\rm if}\ \l_{\rm min}(A)\geq 1 \\
\l_{\rm max}(A) \qquad\ \  {\rm if}\ \l_{\rm max}(A) \leq -1. \\
\end{cases}
$$
In between these two sets the level lines of $\ff$ in $(\l_{\rm min}, \l_{\rm max})$-space
are rays which swing from horizontal to vertical.

An explicit form of this  operator  can be given as $\ff(A) = \vf(\l, \L)$
where $\l=\l_{\rm min}(A)$ and $\L=\L_{\rm max}(A)$, with 
$$
\vf((\l, \L) \ =\ \l \quad {\rm if} \ \ \l\ \geq\ 1
\and
\vf((\l, \L) \ =\ \L  \quad {\rm if} \ \ \L\ \leq\ -1
$$
(as above), and in the region $\l\leq 1, \L \geq -1$ (with $\l\leq \L$) one has
$$
\vf(\l, \L) \ =\ \l \cos \theta + \L \sin \theta
\qquad{\rm with} \quad \cos\theta \ =\ {\L+1  \over \sqrt{(\l-1)^2 + (\L+1)^2}}
$$
\centerline{\qquad 
           \includegraphics[width=.6\textwidth,angle=270,origin=c]{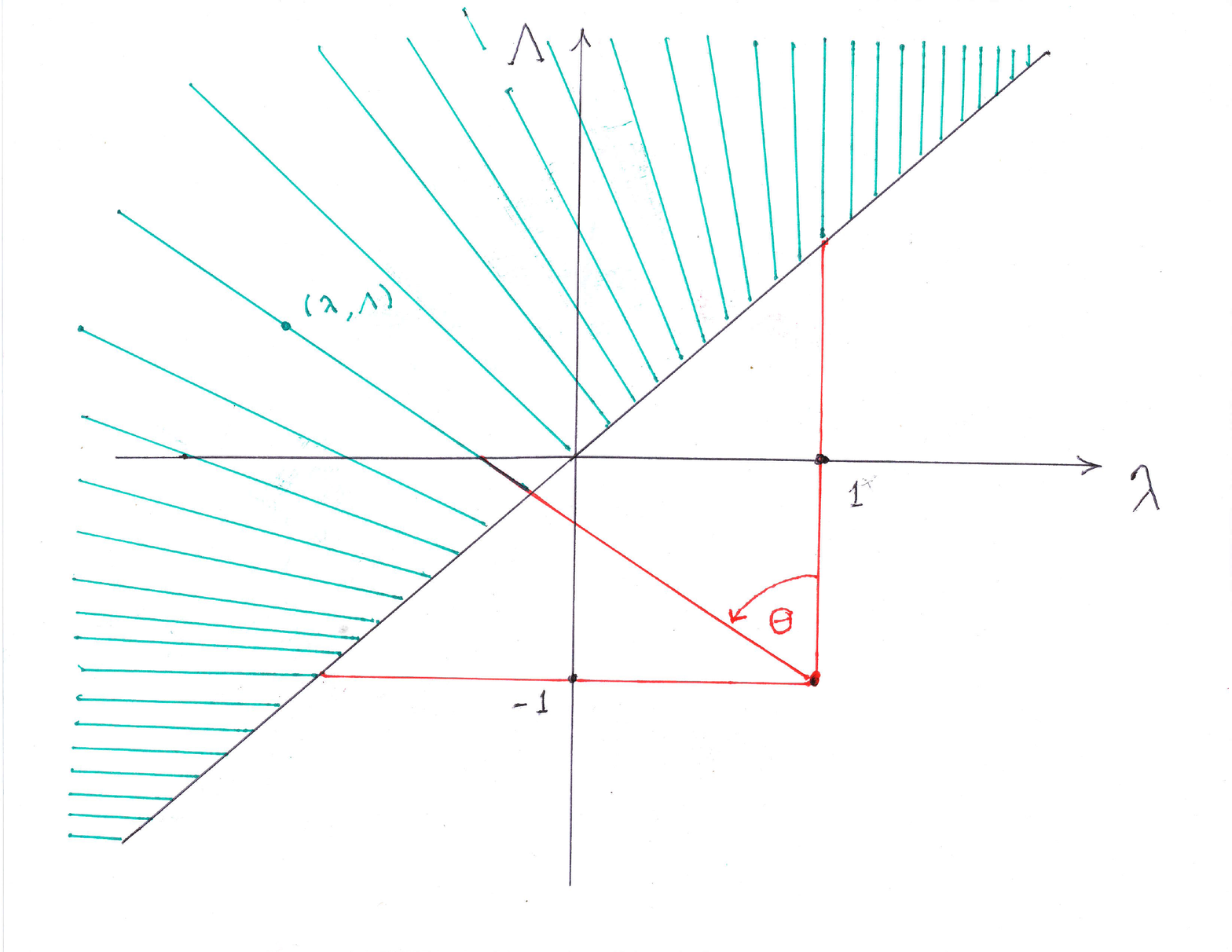}
          }

 Note that this operator  {\bf forces a solution
$u$ to oscillate between being convex and concave} as $\psi$ oscillates between being
$\geq 1$ and $\leq -1$.

\vskip .3in

\centerline{\bf Asymptotic Interiors.}

Let $\ff:\Symn\to \bbr$ be a degenerately elliptic G\aa rding  polynomial of degree $k$ with G\aa rding cone $\G$.
Then for all $c\geq 0$,
$$
\Int\oa \bbf_c \ =\ \G
\and
\overline  \G\ \ss\ \wt{{\overline \G}}
$$
so that a boundary $\bo$  satisfies the $\bbf_c$ (in fact, the $\bbf$) strict boundary hypothesis (Def. \AB.1)
if and only if it is strictly $\bbf = \bbf_0 = \overline \G$-convex.

Let $\G = \G_1  \ss \G_2 \ss \cdots \ss \G_k$ be the(interior) branches of  $\ff$.  Then
 for $\ff \equiv \l_\ell\cdots \l_k$ (G\aa rding eigenvalues), we have  
 $$
 \oa \bbf_c = \overline{\G_\ell}    \and 
 \wt \G_\ell   \ =\ \overline \G_{k -\ell+1}.
$$
Hence, $\bo$ satisfies the strict boundary hypothesis if 
it is strictly $\overline \G_m$-convex for $m=\min\{\ell, k-\ell+1\}$.

Let $\ff:\Symn\to \bbr$ be the canonical operator for a subequation $\bbf$.  Then 
$$
\Int\oa \bbf_c \ =\Int \oa \bbf  \qquad {\rm for \ all \ \ } c\in\bbr.
$$
so strict $\oa  \bbf$-convexity and  strict $\oa {{\wt \bbf}}$-convexity of $\bo$  give the strict boundary hypothesis 
for any inhomogeneous term $\psi$.

\def\ftil{{\tau}}

It is worthwhile to look at computing $\Int \oa \bbf_c$ from $\ff$.   
We have a function $\ftil : \{\Symn - \cp\} \to [-\infty, \infty]$ of degree 0 (i.e.,  a function on the unit sphere in 
$\Symn$), given by
$$
\ftil (A) \ \equiv \liminf_{t\to \infty} \ff(tA) \qquad \text{if $tA \in\bbf$ for all $t\geq $ some $t_0$}
\eqno{(\CC.17)}
$$
and 
$$
\ftil (A) \ \equiv \  -\infty \qquad{\rm otherwise}.
\eqno{(\CC.17)'}
$$

\Lemma {\CC.23}  {\sl  We have
$$
\Int \oa \bbf_c \ \ss\ {\rm cone}\{\ftil > c\}
\eqno{(\CC.18)}
$$
Furthermore, if $\ftil $ is lower semi-continuous, equality holds in (\CC.18).
}

\pf Suppose $A\in \Int \oa \bbf_c$.  Then by Cor.  5.10 of [DD] we have that there exists
$\e>0$ and $R>1$ such that
$$
\text{ $t(A-\e I) \in  \bbf_c$ for all $t\geq R$}
\eqno{(\CC.19)}
$$
Now (\CC.19) means that
$$
\ff(t(A-\e I)) \ \geq\ c \quad \text{for all  $t\geq R$}.
\eqno{(\CC.19)'}
$$
Hence, by the tameness of $\ff$, 
$$
\ff(tA) \ >\ \ff(tA -t\e I) + c(\e) \ \geq\ c  + c(\e) \quad \text{for all  $t\geq R \ (>1)$}.
$$
From (\CC.17) we see that $\ftil (A)>c$, and we have established (1).

Now suppose that $\ftil $ is lower semi-continuous.  Then $\{\ftil >c\}$ is open.
Hence if $\ftil (A)>c$, then $\ftil (A-\e I)>c$ for all $\e>0$ sufficiently small.
This means by (\CC.17) that 
$$
\liminf_{t\to \infty} \ff(t(A-\e I)) \ >\ c,
$$
which by (\CC.19)$'$ (Cor. 5.10 of [\DDD]) means that $A\in \Int  \bbf_c$.
\qed

\Note{\CC.24}  One can rephrase (\CC.18) as
$$
\Int \oa \bbf_c \ \ss\ \bigcup _{\e>0} \oa \bbf_{c_\e}.
$$
Moreover, one can show that 
$$
\{ \ftil >c\} \ =\  \bigcup _{\e>0} \oa \bbf_{c_\e}.
$$
From this one can prove that $\ftil $ is not always l.s.c. (let $\bbf = \cp-I$ and $\ff(A)= \det(A-I)$).
However, it one replaces $tA$ in (\CC.17)  by $tA -\l I$ for $\l>0$ large, lower semi-continuity might be true.

\vskip.3in

\centerline{\bf Further Examples}

\Ex{\CC.25} Consider the operator $\ff$  on $\bbf = \D = \{\tr A\geq0\}$ given as follows.
For $c \geq 1$ the set $$\ff^{-1}(c) =   c I + \partial \cp$$
For $0< c <1$ the set $\ff^{-1}(c) $ has two pieces. We shall use coordinates  $(x,t)\in \{\tr = 0\} \oplus \tr$. 
The first piece is 
$$
\{|x| \leq r(t)\} \times \{t\}
$$
where $r(t)\to\infty$ as $t\to 0$.  The second piece is the part above trace = t of the downward translate
$$
-\rho I + \partial \cp
$$
with $\rho$ chosen so that this set contains the boundary of the ball above.

\ 

\centerline{\qquad 
           \includegraphics[width=.4\textwidth,angle=90,origin=c]{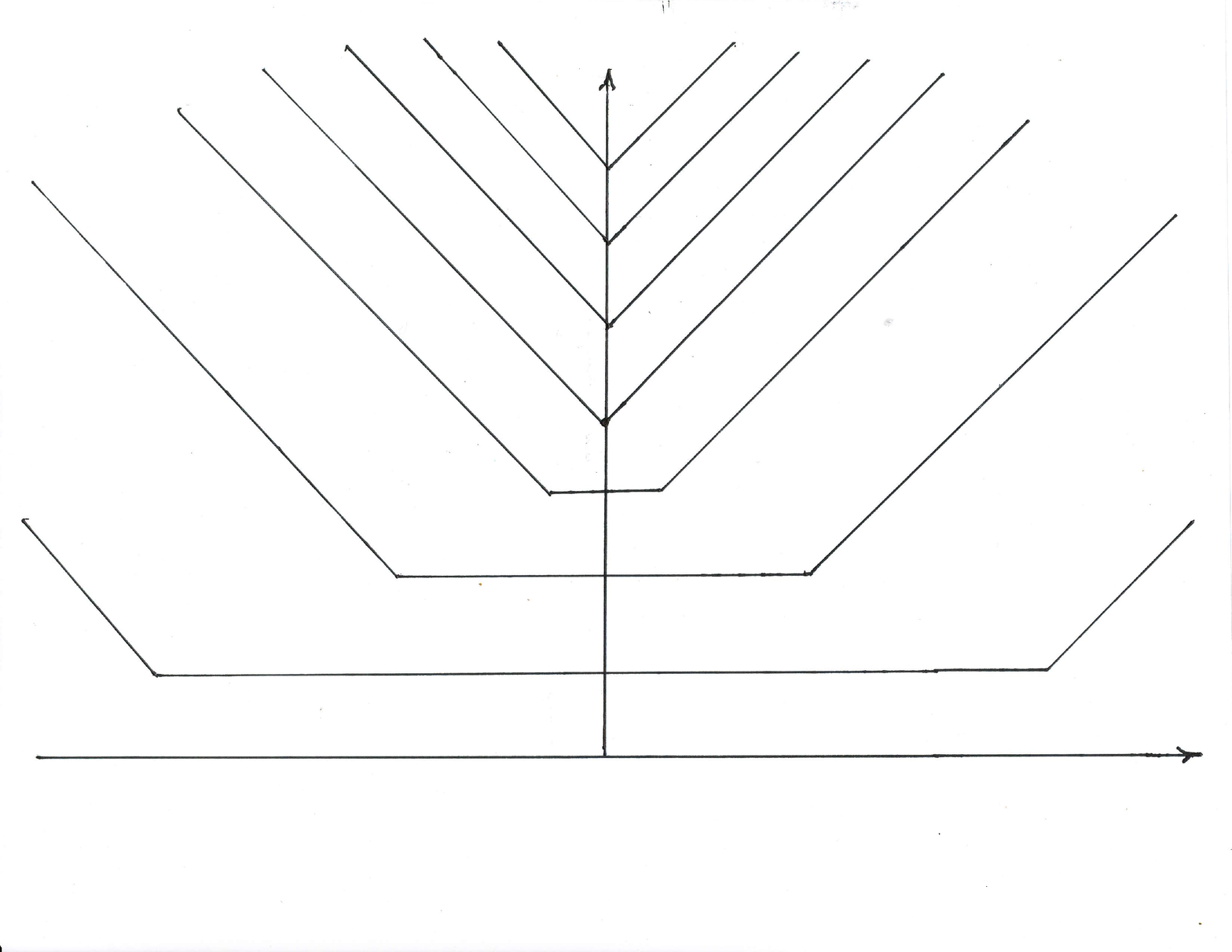}
          }          
Now we have the sets $\bbf_c \equiv \{\ff\geq c\}$, and one computes that 
$$
\Int \oa\bbf_0 \ =\  \Int \D  \and  \Int \oa \bbf_c \ =\ \Int\cp \quad {\rm for}\ \ c>0.
$$

\Ex{\CC.26}  One could expand this by adding the hyperplanes $\{\tr = t\}$ for $-1\leq  t\leq 0$.  Then 
$$
\Int \oa\bbf_c \ =\ \Int \D\quad {\rm for}\ \ -1\leq c \leq 0
$$
One could continue for $t\leq -1$ by inverting what was done for $t\geq 1$, and 
$$
\Int \oa \bbf_c \ =\ \Int \cpt   \quad {\rm for}\ \ c\leq -1.
$$
This can, in fact,  be done for any finite number of jumps.

\Ex{\CC.27} For a continuous example on $\D$, let the part for $c\geq 1$ be as above.
Then between 1 and 0 let the cone open up from $\cp$ to all of $\D$ as $c\downarrow 0$.
Let $\cp(c)$ be the cone with vertex $(0,c)$.  Then
$$
\Int \oa\bbf_0 \ =\ \Int \D, \quad  \Int  \oa \bbf_c \ =\ \Int \cp(c),\ \ 0<c\leq 1 \quad{\rm and}\quad
\Int  \oa \bbf_c\ =\ \Int \cp, \ \ c\geq 1.
$$

\vskip .3in

\noindent{\headfont \CCC.\  Fundamental Solutions.}

It is  natural to ask whether  it is possible  to solve the inhomogeneous Dirichlet problem
$\bdf (D^2 u) = \psi$ where $\psi$ is more general than continuous, for example, a measure.  In 
this section we shall address the basic case where $\psi$ is any (positive) multiple
of the delta function,

We begin with a clear formulation of this problem.  Let $\bbf\ss\Symn$ be a   cone subequation
(with the origin as vertex) which is ST-invariant, i.e., invariant  under a subgroup $G\ss {\rm O}(n)$ which
acts transitively on the sphere $S^{n-1}\ss\rn$.  
Then we fix a degenerate elliptic  operator $\bdf \in C^\infty(\bbf)$ which is $G$-invariant, homogeneous of some degree $m>0$
and $\partial \bbf = \{\bdf =0\}$.
We want to, in some sense, solve the equation 
$$
\bdf (D^2_x K) \ =\ c \d_0 \qquad (c>0) \ \  {\rm on}\ \ \rn.
\eqno{( \CCC.1)}
$$

Now in the situation we are in  (where $\bbf$ is a  ST-invariant cone subequation) there is a natural candidate
for a solution to this problem.  Each such $\bbf$ has attached an invariant {\bf Riesz characteristic} 
$p= p_\bbf \in [1,\infty]$ which is typically easy to compute, and for most interesting subequations
it is finite (see [\TANG \S 3] for discussion and  [\TANG \S 4] for examples).  In fact if $\bbf \equiv \overline \G$ is the closure of the G\aa rding cone $\G$ for a G\aa rding/Dirichlet polynomial $\bdf$, 
then $p\in [1,n]$ since $\bbf\ss\D$ (see  (6.3) and (6.4) in [\TANG]).

Now with $p$ finite, the {\bf Riesz kernel} 
$$
K(x) \ \equiv\ 
\begin{cases}
{1\over 2-p} |x|^{2-p}    \qquad {\rm for} \ \ 1\leq p<2  \\
\quad  \log |x|   \ \ \     \qquad {\rm for} \ \ p=2  \\
- {1\over p-2}{1\over |x|^{p-2}}  \qquad {\rm for} \ \ p>2  \
\end{cases}
\eqno{( \CCC.2)}
$$
is $\bbf$-harmonic in $\rn-\{0\}$ in the viscosity sense  (and $\bbf$-subarmonic across $\{0\}$
since $K$ has no test functions at 0).

Notice that we have not yet mentioned the operator $\bdf$.  For the ordinary inhomogeneous 
Dirichlet problem, with continuous right hand side $\psi\geq 0$,  
we can replace the operator $\bdf$ by any power $\bdf^\a, \a>0$.
That is,  we can replace $\bdf(D^2 u) = \psi$ with 
$\bdf(D^2 u)^\a = \psi^\a$,  and solutions of one are solutions of the other.
However,  for the problem we are now addressing there is one, and only one, exponent $\a$
that solves the problem.

There is a natural way to smooth the Riesz kernel $K$ with a pointwise decreasing family
$K_\e$  of $\bbf$-subharmonics.  Define $k(t)$ so that $k(|x|) = K(x)$ in (\CCC.2).  Set
$$
K_\e(x) \ \equiv\ k\left(  \sqrt{|x|^2+\e^2} \right)
\eqno{( \CCC.3)}
$$
Note that $k(t)$ is increasing for all $p$.  In fact,
$$
k'(t) \ =\ {1\over t^{p-1}}\quad\text{for all $1\leq p <\infty$.}
\eqno{( \CCC.4)}
$$
Hence,
$$
K_\e \in C^\infty(\rn) \ \ \text{decreases pointwise in $\rn$ to  $K$}.
\eqno{( \CCC.5)}
$$

\Lemma { \CCC.1}  
$$
D^2_x K_\e \ =\ {1\over \left ( \sqrt{|x|^2+\e^2}\right )^p}
\left[   P_{x^\perp} - (p-1)P_x + {\e^2 p\over |x|^2+\e^2} P_x    \right]
\ =\ 
{1\over \e^p} D^2_{\left({x\over \e}\right)} K_1.
$$
{\sl and}
$$
D_x K_\e \ =\ {x\over \left(\sqrt{|x|^2+\e^2}     \right)^p}
\ =\ 
{1\over \e^{p-1}} D_{\left({x\over \e}\right)} K_1.
$$

\Cor{ \CCC.2} {\sl The function $K_\e$ is $\bbf$-subharmonic on $\rn$.}

\noindent
{\bf Proof of Corollary  \CCC.2}  By definition of finite Riesz charateristic $p$ we have
$P_{x^\perp} - (p-1)P_x  \in\partial\bbf$ for all $x\neq 0$.  Hence, by degenerate ellipticity
(positivity) of $\bbf$, adding a positive multiple of $P_x$ keeps you in $\bbf$.
Thus $D^2_x K_\e \in\bbf$ for $x\neq 0$, and since $K_\e$ is smooth this also holds at 0.\qed

\noindent
{\bf Proof of Lemma   \CCC.1}  
We use the following formula for the second derivative of a radial function $G(x)=g(|x|)$,
$$
D^2_x G \ =\ {g'(|x|) \over |x|} P_{x^\perp} + g''(|x|)P_x,
\eqno{( \CCC.6)}
$$
applied to $g_\e(t) \equiv k(\sqrt{t^2+\e^2})$.   By ( \CCC.4) we see that
$$
g_\e'(t)\ =\ {t\over (\sqrt{t^2+\e^2})^p}.
\eqno{( \CCC.7)}
$$
Hence, we have
$$
g_\e''(t)\ =\ {1\over (\sqrt{t^2+\e^2})^p}  \left( 1 - {pt^2  \over t^2 +\e^2}     \right).
\eqno{( \CCC.8)}
$$
The formulas for $D^2_x  K_\e$ follow easily  from ( \CCC.6), ( \CCC.7) and ( \CCC.8), and noting that
$$
{1\over \sqrt{|x|^2 + \e^2}}\ =\ {1\over \e} {1\over \sqrt{{|x|^2\over \e^2} + 1} }
\quad{\rm and}\quad
{\e^2\over |x|^2+\e^2} \ =\ {1\over {|x|^2 \over \e^2}+1}.
\ \ \mathqed
$$

\Theorem{ \CCC.4}  {\sl
Suppose $\bbf$ is a conical ST-invariant subequation of finite Riesz characteristic $p$, $1\leq p<\infty$ in $\rn$, and 
let $\bdf\in C^\infty(\bbf)$ be  homogeneous of degree $m>0$ and compatible with $\bbf$. Recall that (\CCC.5) holds.
If we set
$$
\a\equiv {n\over mp} 
\and
\vf(|x|) \ \equiv\ \bdf^\a(D^2_x K_1)
$$
i.e., $p \,\deg( \bdf^\a )=n$, then (and only  then, see ( \CCC.11))
$$
\bdf^\a\left( D^2_x K_\e    \right)  \ =\ {1\over \e^n } \vf\left({|x|\over \e}     \right)
\ \equiv\ \vf_\e(x)
\eqno{( \CCC.9)}
$$
is integrable on $\rn$ and defines a (positive) radial approximate delta function with coefficient $c = \int_{\rn}\vf(|x|)$.
In other words}
$$
``\bdf^\a (D^2K) \ =\ c\d_0".
$$

\pf
By Lemma \CCC.1
$$
\bdf^\a\left(D^2_x K_\e  \right)  \ =\   {1\over \e^{\a p m}}  \bdf^\a\left(D^2_{\left({x\over \e}\right)}  K_1 \right)
\eqno{(\CCC.10)}
$$
which by the definition of $\vf(|x|)$ equals $ {1\over \e^{\a p m}} \vf({|x|\over \e})$.
Since $K_1(x) \in C^\infty(\rn)$ is $\bbf$-subharmonic on $\rn$ by Corollary \CCC.2, and $\bdf\geq 0$
on $\bbf$, we have $\vf(|x|)\geq 0$.

 Lemma  \CCC.5  below states that  $\vf(|x|)$ is integrable on $\rn$, thus completing the proof
 that $``\bdf^\a (D^2K) = c\d_0"$.  Notice that for any value of
$\a$ other than $\a = n/mp$ we have
$$
 \bdf^\a \left(D^2_x K_\e  \right) \ =\ {\e^\d\over \e^n}  \vf\left({|x|\over \e}   \right) \quad\text{with $\d=n-\a mp$},
\eqno{( \CCC.11)}
$$
and the limit of the integral as $\e\to 0$ will be either 0 or $\infty$. Together with Lemma  \CCC.5, this completes the proof of Theorem \CCC.4.\qed

\Lemma {\CCC.5}  {\sl For $\a\equiv n/mp$ one has that}
$$
\vf(|x|) \in C^\infty(\rn) \cap L^1(\rn)
\and
\vf(|x|) > 0.
$$

\pf
By Lemma \CCC.1 with $r\equiv |x|$ and $\e=1$,
$$
\begin{aligned}
\vf(r) \ &\equiv\ 
\bdf^\a   \left( D_x^2 K_1 \right)  \\ 
&=\ 
{1\over (\sqrt{r^2+1} )^n }  \bdf^\a \left(P_{x^\perp} - (p-1) P_x + {p\over r^2+1} P_x       \right).
\end{aligned}
\eqno{(\CCC.12)}
$$
By definition of the Riesz characteristic of $\bbf$, $A\equiv P_{x^\perp} - (p-1) P_x \in \partial \bbf$ and 
$A + {p\over r^2+1} P_x\in \Int \bbf$ for all $0\leq r <\infty$.
Now invoking $(\bbf, \bdf)$ compatibility, we see that $\vf(|x|) >0$.  Since $\bdf$ is $C^\infty$ on $\bbf$,  
 it follows that $\vf(|x|)\in C^\infty(\rn)$.
Set $t\equiv {1\over r}$.  Then $\bdf(A+{p\over r^2+1} P_x) = \bdf(A+{p t^2\over 1 + t^2} P_x)$
is smooth at $t=0$ and equals $\bdf(A) = 0$.
Hence, $\bdf(A+{p\over r^2+1} P_x) \leq C {1\over r}$ for some $C>0$, which proves that
$$
\vf(r) \ \leq \ {C^\a \over r^\a \left(\sqrt{r^2+1}       \right)^n}
$$
so that $\vf(|x|) \in L^1(\rn)$.\qed

\medskip
Now we examine a list of operators $\bdf$, with the powers $\a$ so that  $\bdf^\a (D^2K)=c\d_0$, taken from
\S \CC.

\Ex{ \CCC.6}  In Examples \CC.1 -- \CC.5  where $\bdf$ is the determinant or the $k^{\rm th}$
Hessian operator  (over $\bbr$, $\bbc$ or $\bbh$) the power $\a = 1$.  From this point of view these are very natural operators.

\Ex{ \CCC.7}  In Examples \CC.6 (a), (b) and (c) the operator $\bdf$, equal to the product of the 
$p$-fold sums $\l_I$, is of degree ${n_0\choose p}$ for the cases $\bbr^{n_0}, \bbc^{n_0}$ and  $\bbh^{n_0}$.
One calculates that 
$$
\a\ =\ {1\over {n_0-1\choose p-1}}
$$
is the correct power for $\bdf$.

\Ex{ \CCC.8}
In Example \CC.7  the Lagrangian operator $\bdf$ on $\bbc^{n_0}$
 has degree $m=2^{n_0}$ and Riesz characteristic $p=n_0$
and it should be raised to the power
$$
\a \ =\ {1 \over 2^{n_0-1} }
$$

\Ex{ \CCC.9}
In Example \CC.9  $\bdf_\d$ should be raised to the power 
$$
\a \ =\ {n+\d \over n(1+\d)}.
$$

\Ex{ \CCC.10}
In Example \CC.10 with $\cp_{\l, \L}^-$ the Pucci extremal operator,  this operator should be raised to the power
$$
\a \ =\ {\l \over \L} (n-1) -1
$$

For the  G\aa rding-Pucci operator  $\bdf_{\l,\L}$ one has $p_\bbf = {\l\over \L}(n-1)+1$, $m=\deg \bdf_{\l,\L} = 1$, and so 
$\a = n( {\l\over \L}(n-1)+1)^{-1}$

\vskip.3in



\noindent{\headfont \DD. \ Two Important Questions.}

\def\bfpsi{{\bbf_\bdf(\psi)}}

We pose two questions concerning topologically tame operators $\bdf\in C(\bbf)$
with constant coefficients in $\rn$.

\noindent
{\bf Question A.}  Does comparison always hold for $\bfpsi$?

We note that comparison holds  for $\bfpsi$  if and only if the following
subaffine  property holds:
$$\left\{
\begin{aligned}
& \text{$u$ is $\bfpsi$-subharmonic and } \\
&\text{  $v$ is $\wt{\bfpsi}$-subharmonic}
\end{aligned}\right\}
 \quad\Rightarrow \quad 
 \text{$u+v$ is subaffine.}
$$
More specifically, does comparison hold for $\tr\{\arctan D^2u\}= \psi$ (Example \CC.18)?
Any counterexample cannot be tamable, so Examples \CC.20 and \CC.21 also provide
candidates.  Also the functions $u$ and $v$ cannot be quasi-convex since if they are, then the subaffine property holds.
(Note that Lemma A.2 on quasi-convex approximation  requires tameness.)

The second question can be stated succinctly as follows: 
{\sl Is what we are calling  an operator $\bdf$ actually single-valued if it is topologically tame,
or even if it is tame?}

More precisely, given an open set $X\ss\rn$, define the {\sl functional domain} of $\bdf$ to be 

\smallskip
\centerline
{
$\cd_X(\bdf) \ =\ \{u : u\ \text{is $\bfpsi$-harmonic for some $\psi\in C(X)\}$}$
}
\noindent
and define $\bdf(D^2 u)$, for $u\in \cd_X(\bdf)$ as the set of $\psi\in C(X)$ such that 
$u$ is $\bfpsi$-harmonic.  Then

\noindent
{\bf Question B.} Is $\bdf$ single-valued on $\cd_X(\bdf)$?

If $\bdf$ is uniformly elliptic and convex, then it is well known that $\bdf$ is single-valued.
Nothing more seems to be known.  We have asked several experts this question.
The little more we can say is the following.

First the general case, with $\psi\in C(X)$ arbitrary, is equivalent to the case
where $\psi$ is constant.  Hence the question can be restated as

\noindent
{\bf Question B$'$.}   Can a pair of subequations ${\bf H}, \bbf \ss\Symn$ with
${\bf H} \ss \Int \bbf $ have a simultaneous harmonic $u$?

Since  $\partial {\bf H}$ and  $\partial  \bbf$ are disjoint, a counterexample must fail
to be quasi-convex in all neighborhoods of all points by quasi-convex addition.

\Prop{\DD.1} {\sl Each canonical operator $\bdf$ (Proposition \CC.13)
  is single-valued.}


\vskip .3in

\noindent{\headfont  Appendix A.\  Comparison for Constant Coefficient Operators}

\def\BB{A}

The uniqueness part of Theorem \AA.7 holds for any domain $\O\ss\ss\rn$ without the 
assumption of boundary convexity.   The argument for this comparison result is easier than the one given 
for the Main Theorem \AA.11, and so we are including it here.

\Theorem{\BB.1}  {\sl
Let  $(\bbf, \bdf)$ be a reduced  subequation-operator pair.
Suppose that   the operator $\bdf\in C(\bbf)$ can be tamed and that $\psi\in C(\ob)$
takes values in $\bdf(\bbf)$ (i.e., is admissible).
Suppose  $u, v\in\USC(\ob)$ with $u$  $\bfpsi$-subharmonic 
and $v$  $\wt \bfpsi$-subharmonic on $\O$. Then}
$$
{\rm If}\ \ u+v \ \leq\ 0 \ \ {\rm  on}\ \ \bo, \ \ {\rm then}\ \ \ \ u+v \ \leq\ 0 \ \ {\rm  on}\ \ \ob.
\eqno{(\BB.1)}
$$

\medskip

\def\FA{{\mathbb F}}

We shall give two proofs.  The first is based on the Theorem on Sums of Crandall, Iishi and
Lions [\CIL], which in turn is based on the Slodkowski/Jensen Lemma.  The second is based on
an Almost Everywhere Theorem and the notion of a subaffine function.  This A.E.\ Theorem
also rests on the same Slodkowski/Jensen Lemma (see [\AET]).  These proofs provided the 
original motivation for the concept of tameness.  
Without the tameness of the operator, the old arguments did not apply.

\medskip
\noindent
{\bf Proof I. Step 1 (Strict Approximation)}. 
For this first  proof we  simplify the notation for $\bfpsi$ to $\FA$, suppressing the 
dependence on both $f$ and $\psi$.

Consider
$$
u_\l(x) \ \equiv \ u(x) + {\l\over2} |x|^2
\eqno{(\BB.2)}
$$
for $\l>0$.  Note that $\vf$ is a test function for $u$ at $x_0$
 $\iff$ $\vf_\l(x)  \equiv  \vf(x) + {\l\over2} |x|^2$ is a test function
 for $u_\l$ at $x_0$, and
 $$
D^2_{x_0} \vf_\l =  D^2_{x_0} \vf + \l I.
\eqno{(\BB.3)}
$$

Define $\FA^\l$ by its fibres
$$
\FA^\l_x \ \equiv \ \FA_x  + \l I.
\eqno{(\BB.4)}
$$
This $\FA^\l$ is a subequation, and we can restate  (\BB.3) by saying
$$
\text{$u$ is $\FA$-subharm. $\iff$
  $u_\l = u+{\l\over 2}|x|^2$ is
$\FA^\l$-subharm.}
\eqno{(\BB.3)'}
$$

Note that a function  $u$ is $\FA$-subharmonic 
if and only if $u+c$ is $\FA$-subharmonic for all $c\in \bbr$ since $\FA$ is reduced.
Thus we may assume that ``0'' in (\BB.1) can be replaces by any constant $c$.
Now,  since $u_\l$ decreases to $u$ as $\l\downarrow 0$,
 it suffices to prove the theorem with $u$ replaced by $u_\l$. 
That is, we assume that $u$ is $\FA^\l$-subharmonic
for some $\l>0$.

\medskip
\noindent
{\bf Step 2 (Calculating the Dual)}.  From (\AA.3) we see that the fibres of the dual subequation are given by
$$
\wt \FA_y \ = \ \wt \bbf \cup \{B : -B \in \Int \bbf \ {\rm and}\ \bdf(-B) \leq \psi(y)\}.
\eqno{(\BB.5)}
$$

The final step is the main step.

\medskip
\noindent
{\bf Step 3. (Apply the Theorem on Sums [\CIL]).}
The statement we draw on is the following, given in Theorem C.1 in [\DDR].
If $u+v$ has an interior maximum at $x_0\in\O$ which is strictly larger than
the maximum on $\bo$, then there exist: 

\noindent
(1) \ \  numbers
$\e\downarrow 0$ and points $(x_\e, y_\e)\in\O\times\O$ such that  

\centerline{ $(x_\e, y_\e) \ \to\ (x_0,x_0)$ as $\e\downarrow 0$,}

\noindent
(2) \ \ $A_\e \in \FA^\l_{x_\e}$,\quad and \qquad  (3) \ \ $B_\e \in\wt { \FA}_{y_\e}$,

 such that

\noindent
(4)\ \  $A_\e+B_\e \ \leq\ 0.$

Set $P_\e \equiv -(A_\e+B_\e)$ so that we can replace (4) by

\noindent
(4$'$) \ \ $-B_\e \equiv A_\e+P_\e$ with $P_\e\geq0$. 

Now by the definition of $\FA^\l$ and positivity, condition (2) states that

\noindent
(2)$'$ \ \ $A_\e+P_\e - \l I  \in\bbf$ \ \ and\ \ $\bdf(A_\e+P_\e - \l I ) \ \geq\ \psi(x_\e)$.

By (\BB.5) condition (3) states that
$$ \begin{aligned} \text{either}\  &\text {(a)\ \  $A_\e+P_\e =-B_\e \notin \Int \bbf$} \\
\text{ or}\quad &\text{(b) \ \ $A_\e+P_\e = -B_\e \in \Int \bbf$ \ and \ 
$\bdf(A_\e+P_\e)\leq \psi(y_e)$.  }
\end{aligned}
\leqno{(3')}$$

\noindent
Now $\bbf +\l I \ss\Int \bbf$ so that (3a)$'$ is ruled out by (2)$'$.  Thus, the inequality in  (3b)$'$ must hold.  
With 
$$
A_\e' \  \equiv \ A_\e +P_\e -\l I
$$
we now see that the combination of conditions (2), (3) and (4) (or equivalently (2)$'$, (3)$'$ and (4)$'$)
are equivalent to the single condition:

\noindent
(5) \ \ $A_\e' \in\bbf$ \ \ \ and\ \ \   $\psi(x_\e)\  \leq \ \bdf(A_\e') \  \leq \  \bdf(A_\e'+\l I) \ \leq \ \psi(y_\e)$.

Taking the limit as $\e\downarrow 0$, we see that the tameness assumption on the operator $\bdf$ yields the contradiction.\qed

\medskip
\noindent
{\bf Proof II. (An Outline).}  Some readers may find this proof to have clearer motivation and
more intuitive appeal.  In addition, this proof establishes quasi-convex approximation for the 
subequations $\bfpsi$ and $\wt\bfpsi$ even though they do not have constant coefficients.

\noindent
{\bf Step I.}  Show that if $u$ and $v$ are $C^2$, then $D^2u+D^2v \in \cpt$.  That is,
$w\equiv u+v$ is $\cpt$-subharmonic where $\cpt \equiv \{A  : \l_{\rm max}(A)\geq0\}$
is the dual of the subequation $\cp\equiv  \{A  : \l_{\rm min}(A)\geq0\}$.
This is an algebraic step which is valid in much greater generality. Namely, for any
closed subset ${\bf G} \ss\Symn$,
$$
\text{if ${\bf G} +\cp \ \ss\ {\bf G}$, then ${\bf G} + \wt {\bf G} \ \ss\ \cpt$.}
\eqno{(\BB.6)}
$$

\medskip
\noindent
{\bf Step II.}
Recall from [\DDD] that for an upper semi-continuous
function $w$
$$
\text{$w$ is subaffine $\qquad\iff\qquad w$ is $\cpt$-subharmonic.}
\eqno{(\BB.7)}
$$
Thus the concept of being  ``sub'' the affine functions has an advantage over satisfying the maximum principle (i.e., being  ``sub'' the constants).
It is a local concept.

\medskip
\noindent
{\bf Step III.}  
Suppose $u$ and $v$ are quasi-convex.  Then by Alexandrov's Theorem 
 both are twice differentiable almost everywhere, and we have 
 $D^2_x u\in\bfpsi_x$ and $D^2_x v\in \wt\bfpsi_x$ for almost all $x$.  Therefore by (\BB.5)
 $$
D^2_x(u+v) \ =\ D^2_xu + D^2_x v \ \in \ \cpt \ \ \ \text{for a.a. $x$.}
\eqno{(\BB.8)}
$$

\medskip
\noindent
{\bf Step IV. (Apply the AE Theorem).}  
This result (see [\AET])  states that for any subequation ${\bf G}$ and any
locally quasi-convex function $w$:  \medskip

\centerline
{
If the 2-jet $J^2_x w \in {\bf G}_x$ for a.a. $x$, then $w$ is ${\bf G}$-subharmonic.
}

\medskip
\noindent
{\bf Step V.}  
At this point we have proved the theorem for $u$ and $v$ quasi-convex, so that it suffices
to establish quasi-convex approximation for $\bfpsi$ and $\wt\bfpsi$.

\Lemma {\BB.2. (Quasi-Convex Approximation)}  {\sl If the operator $f$ is tame, and 
if there exist an $\bfpsi$-subharmonic function and an $\wt\bfpsi$-subharmonic function which are bounded 
below, then

(a) Each $\bfpsi$-subharmonic function $u$ can be approximated by a decreasing 
sequence of quasi-convex $\bfpsi$-subharmonic functions $\{u_j\}$  converging pointwise to $u$.

(b) Each $\wt\bfpsi$-subharmonic function $v$ can be approximated by a decreasing 
sequence of quasi-convex $\wt\bfpsi$-subharmonic functions $\{v_j\}$  converging pointwise to $v$.
}

\medskip
\noindent
{\bf Proof of (a).} By replacing $u$ by max$\{u,  \a -N\}$,
where $\a$ is an  $\bfpsi$-subharmonic function which is bounded below,
 we can assume that $u$ is bounded by $M$.
Let  $(u_\l)^\e$ be the strict approximation $u_\l$ in (\BB.2) followed by
the standard $\e$-sup-convolution.
It suffices to show that:
$$
(u_\l)^\e \ \ \text{ is $\bfpsi$-subharmonic if $\e$ is small relative to} \ \l
\eqno{(\BB.9)}
$$

The function $(u_\l)^\e$ is the supremum taken over $|z|\leq\d\equiv\sqrt{2\e M}$ of the 
functions
$$
v(x) \ \equiv\ u(x-z) + {\l\over 2} |x-z|^2 - {1\over \e}|z|^2.
$$
First we show that each $v(x)$ is $\bfpsi$-subharmonic. Suppose that $\vf$ is a 
test function for  $v$ at a point $x_0$.  We must show that  $B\equiv D^2_{x_0}\vf \in \bfpsi_{x_0}$,
i.e., $B\in \bbf$ and $\bdf(B)\geq \psi(x_0)$. Since
$$
v(x)\ \leq\ \vf(x) \ \ \ \text{near \ $x_0$\  with equality at \ $x_0$},
\eqno{(\BB.10)}
$$
if we set $y\equiv x-z$ and $y_0=x_0-z$, it follows that 
$\overline \vf(y) \equiv \vf(y+z)-{\l\over 2}|y|^2 +{1\over \e}|z|^2$ is a test function
for $u(y)$ at $y_0$.  That is, 
$$
u(x)\ \leq\  \overline\vf(x) \ \ \ \text{near \ $y_0$\  with equality at \ $y_0$}.
\eqno{(\BB.10)'}
$$
Hence, $A\equiv D^2_{y_0}\overline\vf \in \bbf$ and $\bdf(A)\geq \psi(y_0)$.
Therefore, $B=A +\l I \in \bbf$.
Let $\o(\d)$ denote the modulus of continuity of $\psi$.
Since $\bdf$ is degenerate elliptic  on $\bbf$, we have
$$
\bdf(B) - \psi(x_0) \ =\ \bdf(A+\l I) -\psi(x_0) \geq \ c(\l) + \psi(y_0) -\psi(x_0)
\ \geq\ c(\l) - \o(\d)
$$
since $|y_0-x_0| = |z| \leq \d \equiv \sqrt{2\e M}$.  With $\e$ small, 
$c(\l)-\o( \sqrt{2\e M})\geq0$ which proves that each $v$ is $\bfpsi$-subharmonic.

The rest of the proof is standard and goes as in the constant coefficient case
(see [\CIL], [\CRA] or [\DDD]).  The proof of (b) is similar.\qed





\vskip .3in

\def\AAA{B}

\noindent{\headfont  Appendix \AAA.\ Compatibility and Topological  Tameness are Necessary Conditions }

\Prop{\AAA.1}   {\sl
Suppose $\bdf$ is an elliptic operator on $\bbf$.    Assuming comparison and  existence for the 
Dirichlet problem $\bdf(D^2 u)=c$ (constant),  implies topological tameness and $\bbff$-compatibility.
}

\pf
In Corollary \CC.17 we proved that comparison for the inhomogeneous equation  
for $\psi \equiv c$  (admissible), implies that $\bdf$ is topologically tame.

We must show that if, in addition, existence holds, then $\bdf$ and $\bbf$ must be compatible.
Suppose they are not compatible.  Then we can assume there exist
$$
c \in \bdf(\partial \bbf), \ A \in \partial \bbf \and \bdf(A) >c.
$$
Set $u(x) \equiv \half \bra {Ax}x$ and $\vf=u\bigr|_{\bo}$.
Suppose existence holds, and fix a small ball $\O\ss\rn$. Let $v$ denote the solution to 
$$
f(D^2 v) \ = \ c, \ \   D^2 v \in \bbf \qquad {\rm with}\qquad v\bigr|_{\bo} \ =\ \vf.
\eqno{(\II.1)'}
$$
By comparison $v$ must equal the  Perron function.  However, the Perron family
for this is the same as the Perron family for 
$$
D^2 u \in \partial \bbf_c, \qquad {\rm with}\qquad u\bigr|_{\bo} \ =\ \vf.
\eqno{(\II.1)}
$$
By Theorem \AA.4 the solution to (\II.1) is the $\bbf_c$-harmonic
$u(x) \equiv \half \bra{Ax}x$ (since $A\in \partial \bbf_c$).
However, this Perron function is not a solution to (\II.1)$'$ since
$u$ is a $C^2$-function with 
 $\bdf(D^2 u)=\bdf(A)>c$. \qed

\Note{\AAA.2}  For  the 
Dirichlet problem $\bdf(D^2 u)=c$ on $\bbf$-convex domains in $\rn$,   
topological tameness and $\bbff$-compatibility are also sufficient.


\def\item{}
\vskip .3in

\centerline{ \headfont References}

\vskip .2in

\noindent
\item{[\ALL]}  S. Alesker,   {\sl  Quaternionic Monge-Amp\`ere equations}, 
J. Geom. Anal., {\bf 13} (2003),  205-238.
 ArXiv:math.CV/0208805.

\noindent
\item{[\AV]} 
 \ \----------  and M. Verbitsky,  {\sl Quaternionic Monge-Ampre equation and Calabi problem for HKT-manifolds}, Israel J. Math., {\bf 176} (2010), 109Ð138.

\smallskip

\noindent
\item{[\BTA]}   E. Bedford and B. A. Taylor,  {The Dirichlet problem for a complex Monge-Amp\`ere equation}, 
Inventiones Math.{\bf 37} (1976), no.1, 1-44.

\smallskip

\noindent
\item{[\BL]}
Z. Blocki, {\sl
Weak solutions to the complex Hessian equation}, 
Ann. Inst. Fourier (Grenoble) {\bf 55}  (2005), no. 5, 1735-1756.

\smallskip

\noindent
 \item{[\BRE]}  H. J. Bremermann,
    {\sl  On a generalized Dirichlet problem for plurisubharmonic functions and pseudo-convex domains},
          Trans. A. M. S.  {\bf 91}  (1959), 246-276.
\smallskip

\noindent
 \item{[\CNS]}  L. Caffarelli, L. Nirenberg and J. Spruck,  {\sl
The Dirichlet problem for nonlinear second order elliptic equations, III: 
Functions of the eigenvalues of the Hessian},  Acta Math.
  {\bf 155} (1985),   261-301.

 \smallskip

\noindent
\item{[\CP]}  M.  Cirant and K. Payne, {\sl On viscosity solutions to the Dirichlet problem for elliptic branches of 
inhomogeneous fully nonlinear equations}, Publ.  Mat. {\bf 61} (2017), 529-575.
 \smallskip

\noindent
\item{[\CPW]}  T. C. Collins, S. Picard and S. Wu, {\sl Concavity of the Lagrangian phase operator
and applications}, ArXiv:1607.07194v1.
 \smallskip

\noindent
\item{[\CIL]}   M. G. Crandall, H. Ishii and P. L. Lions {\sl
User's guide to viscosity solutions of second order partial differential equations},  
Bull. Amer. Math. Soc. (N. S.) {\bf 27} (1992), 1-67.

 \smallskip

\noindent
\item{[\CRA]}   M. G. Crandall,  {\sl  Viscosity solutions: a primer},  
pp. 1-43 in ``Viscosity Solutions and Applications''  Ed.'s Dolcetta and Lions, 
SLNM {\bf 1660}, Springer Press, New York, 1997.

 \smallskip

\noindent
\item{[\DTT]}   S. Dinew, H.-S.\ Do and T.\ D.\ T\^o, {\sl A viscosity approach to the 
Dirichlet problem  for  degenerate complex Hessian type equations}, ArXiv:1712.08572.

 \smallskip

\noindent
\item{[\DON]}   S. Donaldson,  {\sl  Moment maps and diffeomorphisms},  
Asian J. Math. {\bf 3}, No. 1 (1999), 1-16.

 \smallskip

 \noindent
\item{[\CG]}   F. R. Harvey and H. B. Lawson, Jr., {\sl Calibrated geometries}, Acta Mathematica 
{\bf 148} (1982), 47-157.

\smallskip

 \noindent
\item{[\DDD]}   \ \----------,  {\sl  Dirichlet duality and the non-linear Dirichlet problem},    Comm. on Pure and Applied Math. {\bf 62} (2009), 396-443. ArXiv:math.0710.3991.

\smallskip

 \noindent
\item{[\DDR]}  \ \----------, {\sl Dirichlet Duality and the Nonlinear Dirichlet Problem on Riemannian Manifolds},  J. Diff. Geom. {\bf 88} (2011), 395-482.   ArXiv:0912.5220.

\smallskip

 \noindent
\item {[\GEO]}  \ \----------,  {\sl  Geometric plurisubharmonicity and convexity - an introduction},
  Advances in Math.  {\bf 230} (2012), 2428-2456.   ArXiv:1111.3875.

\smallskip

 \noindent
\item {[\HYP]} \ \----------, {\sl  Hyperbolic polynomials and the Dirichlet problem},   ArXiv:0912.5220.

\smallskip

 \noindent
\item {[\BEL]}   \ \----------,  {\sl The equivalence of  viscosity and distributional
subsolutions for convex subequations -- the strong Bellman principle},  
 Bulletin Brazilian Math. Soc. {\bf 44} No. 4 (2013),  621-652.  ArXiv:1301.4914.

\smallskip

 \noindent
\item {[\SURVEY]}  \ \----------,  {\sl  Existence, uniqueness and removable singularities
for nonlinear partial differential equations in geometry},
pp. 102-156 in ``Surveys in Differential Geometry 2013'', vol. 18,  
H.-D. Cao and S.-T. Yau eds., International Press, Somerville, MA, 2013.
ArXiv:1303.1117.

\smallskip

 \noindent
\item{[\AC]}  \ \----------,  {\sl  Potential theory on almost complex manifolds},  Ann.\ Inst.\ Fourier, Vol. 65 no. 1 (2015), p. 171-210. ArXiv:1107.2584.

\smallskip

\noindent
\item  {[\NOTES]}  \ \----------,  {\sl Notes on the differentiation of quasi-convex functions,} 
 \ ArXiv:1309:1772.

\smallskip

\noindent
\item  {[\AET]}   \ \----------, {\sl The AE Theorem and Addition Theorems for quasi-convex functions,} 
 ArXiv:1309:1770.

\smallskip

 \noindent
\item {[\TANG]}  \ \----------,  {\sl  Tangents to subsolutions -- existence and uniqueness, Part I}, 
  to appear in Ann. Fac.de Sciences de Toulouse.   ArXiv:1408.5797.

\smallskip

 \noindent
\item {[\TANGG]}  \ \----------,  {\sl  Tangents to subsolutions -- existence and uniqueness, Part II}, 
 J. Geom. Analysis, {\bf 27} (2017), 2190-2223.  ArXiv:1408.5851.

\smallskip

 \noindent
\item {[\LAG]}   \ \----------, {\sl  Lagrangian potential theory and a Lagrangian equation of Monge-Amp\`ere type},\ 
ArXiv:1712.03525.

\smallskip

 \noindent
\item {[\SLE]}   \ \----------, {\sl  The special Lagrangian potential equation},\ 
(to appear).

\smallskip

   \noindent
\item{[\JTY]}   F. Jiang, N. S. Trudinger, and X.-P. Yang, {\sl On the Dirichlet problem for Monge-Amp\`ere type equations},  Calc. Var. Partial Differential Equations {\bf 49} (2014), no. 3-4, 1223Ð1236.

\smallskip

   \noindent
\item{[\KRY]}    N. V. Krylov,    {\sl  On the general notion of fully nonlinear second-order elliptic equations},    Trans. Amer. Math. Soc. (3)
 {\bf  347}  (1979), 30-34.

\smallskip

 \noindent
\item {[\PLI]}  S. Pli\'s, {\sl The Monge-Amp\`ere equation on almost complex manifolds}, 
ArXiv:1106.3356, June, 2011.

\smallskip

\noindent
\item {[\RT]} J. B. Rauch and B. A. Taylor, {\sl  The Dirichlet problem for the 
multidimensional Monge-Amp\`ere equation},
Rocky Mountain J. Math {\bf 7}    (1977), 345-364.

\smallskip

\noindent
\item {[\SLO]}  Z. Slodkowski, {\sl  The Bremermann-Dirichlet problem for $q$-plurisubharmonic functions},
Ann. Scuola Norm. Sup. Pisa Cl. Sci. (4)  {\bf 11}    (1984),  303-326.

\smallskip

\noindent
\item {[\SPR]}  J. Spruck, {\sl  Geometric aspects of the theory of fully nonlinear elliptic equations},  Global theory of minimal surfaces, 283Ð309,
Clay Math. Proc., 2, Amer. Math. Soc., Providence, RI, 2005.

\smallskip

\noindent
\item{[\TRUU]} 
N. S. Trudinger,  {\sl On the Dirichlet problem for Hessian equations},  Acta Math.  {\bf 175} (1995), 151-164.

\smallskip

\noindent
\item{[\TRU]} 
  \ \----------,  {\sl Weak solutions of hessian equations}, Communications in Partial Differential
Equations, {\bf 22(7\&8)} (1997), 1251-1261. 

\smallskip

 \noindent
\item{[\TWC]}  N. S. Trudinger and X.-J. Wang,  {\sl Hessian Measures I}, Topol. Methods Nonlinear Anal. {\bf 19} (1997), 225Ð239.

\smallskip

 \noindent
\item{[\TWCC]}    \ \----------,  {\sl Hessian Measures II}, Ann. Math. {\bf 150} (1999), 579Ð604.

\smallskip

 \noindent
\item{[\TWCCC]}    \ \----------,   {\sl Hessian Measures III}, Journal of Functional Analysis {\bf 193} (2002), 1Ð23.

\smallskip

\noindent
\item {[\WAL]}   J.  B. Walsh,  {\sl Continuity of envelopes of plurisubharmonic functions},
 J. Math. Mech. 
{\bf 18}  (1968-69),   143-148.

 \smallskip

\end{document}